\documentclass[a4paper, inner=3cm, outer=3cm, top=3cm, bottom=3cm, 12pt]{article}

\usepackage{latexsym,amssymb,amsmath,epsfig,amsfonts,graphicx,mathrsfs}
\usepackage{bm}
\usepackage{pstricks,pst-node,pst-tree}
\usepackage{tikz}
\usetikzlibrary{arrows}

\usepackage{booktabs}

\newcommand{\vc}[1]{\boldsymbol{#1}}

\usepackage{dsfont}
\newcommand{\qed}{\hfill $\square$}

\newcommand{\registered}
   {{\scriptsize \ooalign{\hfil\raise0.07ex\hbox{\scriptsize \sc r}\hfil%
              \crcr\mathhexbox20D}}}

\newtheorem{theorem}{Theorem}
\newtheorem{assumption}{Assumption}
\newtheorem{lemma}{Lemma}

\newtheorem{corollary}{Corollary}

\newtheorem{remark}{Remark}
\newtheorem{proposition}{Proposition}

\newcommand{\nc}{\newcommand}
\nc{\ds}{\displaystyle}
\nc{\mbZ}{\mathbb Z}
\nc{\mbQ}{\mathbb Q}
\nc{\mbR}{\mathbb R}
\nc{\mbC}{\mathbb C}
\nc{\mbN}{\mathbb N}
\nc{\mbE}{\mathbb E}
\nc{\mbP}{\mathbb P}

\nc{\PH}{\emph{PH} }
\nc{\ME}{\emph{ME} }
\nc{\LST}{\emph{LST} }
\nc{\rank}{\mbox{rank\hspace{1pt}}}

\numberwithin{equation}{section}

\numberwithin{figure}{section}

\usepackage[english]{babel}

\begin{document}
\title{Extinction in lower Hessenberg branching processes with countably many types}
\author{P. Braunsteins and S. Hautphenne}
\date{\today}
\maketitle

\begin{abstract}
We consider a class of branching processes with countably many types which we refer to as \emph{Lower Hessenberg branching processes}. These are multitype Galton-Watson processes with typeset $\mathcal{X}=\{0,1,2,\dots\}$, in which individuals of type $i$ may give birth to offspring of type $j\leq i+1$ only.  For this class of processes, we study the set $S$ of fixed points of the progeny generating function. In particular, we highlight the existence of a continuum of fixed points whose minimum is the global extinction probability vector $\bm{q}$ and whose maximum is the partial extinction probability vector $\bm{\tilde{q}}$.
In the case where $\bm{\tilde{q}}=\bm{1}$, we derive a global extinction criterion which holds under second moment conditions, and when $\bm{\tilde{q}}<\bm{1}$ we develop necessary and sufficient conditions for $\bm{q}=\bm{\tilde{q}}$. 
\\\noindent \textbf{Keywords}:  infinite-type branching process; extinction probability; extinction criterion; fixed point; varying environment.
\end{abstract}

\section{Introduction}\label{intro}

Multitype Galton-Watson branching processes (MGWBPs) describe the evolution of a population of independent individuals who live for a single generation and, at death, randomly give birth to offspring that may be of various types.
Classical reference books on MGWBPs include Harris \cite{Har63}, Mode \cite{mode}, Athreya and Ney \cite{athreya_ney}, and Jagers \cite{jagers}.
MGWBPs have been used to model populations in several fields, including in molecular biology, ecology, epidemiology, and evolutionary theory, as well as in particle physics, chemistry, and computer science. 
Recent books with a special emphasis on applications are Axelrod and Kimmel \cite{kimmel}, and Haccou, Jagers and Vatutin \cite{Hac05}. 
Branching processes with an infinite number of types have been used to model the dynamics of escape mutants \cite{Serra07} and the spread of parasites through a host population \cite{Bar94,Bar93}; see also \cite[Chapter~7]{kimmel} for other biological applications of infinite-type branching processes.

One of the main quantities of interest in a MGWBP is the probability that the population eventually becomes empty or \emph{extinct}. 
Let the vector $\vc Z_n=(Z_{n,\ell})_{\ell\in\mathcal{X}}$ record the number of type-$\ell$ individuals alive in generation $n$ of a population whose members take types that belong to the countable set $\mathcal{X}$.
We let 
\begin{equation}\label{global_def}
q_i=\mbP [\lim_{n\to\infty} \textstyle{\sum _{\ell\in\mathcal{X}} }\,Z_{n,\ell}=0\,|\, \varphi_0=i]
\end{equation}
be the probability of \emph{global extinction} given that the population begins with a single individual of type $\varphi_0=i$, and we refer to $\vc q:=(q_i)_{i \in \mathcal{X}}$ as the \emph{global extinction probability vector}.
When the set $\mathcal{X}$ contains only finitely many types, many of the fundamental questions concerning $\bm{q}$ have been resolved. In particular, it is well known that \emph{(i)} $\bm{q}$ is the minimal non-negative solution of the fixed point equation $\vc s=\vc G(\vc s)$,
where $\vc G(\vc s):=(G_i(\vc s))_{i \in \mathcal{X}}$, defined in \eqref{Gs}, records the probability generating function associated with the reproduction law of each type, and that \emph{(ii)} if the process is irreducible, then the set of fixed point solutions
\begin{equation}\label{S}
S
= \{ \vc{s} \in [0,1]^{\mathcal{X}} : \vc{s}=\vc{G}(\vc{s}) \}
\end{equation} 
contains at most two elements, $\bm{q}$ and $\bm{1}$. In addition, there is a well-established extinction criterion, namely $\bm{q}=\bm{1}$ if and only if the Perron-Frobenius eigenvalue of the mean progeny matrix (defined in \eqref{MPM}) is less than or equal to one.

If we allow $\mathcal{X}$ to contain countably infinitely many types then this complicates matters considerably. 
Indeed, even the definition of extinction is no longer unambiguous. 
We let 
\begin{equation}\label{partial_def}
\tilde{q}_i=\mbP[\lim_{n\to\infty} Z_{n,\ell}=0,\;\forall \ell\in \mathcal{X}\,|\, \varphi_0=i],
\end{equation}
be the probability of \emph{partial extinction} given that the population begins with a single individual of type $i$, and we refer to $\bm{\tilde q}=(\tilde q_i)_{i \in \mathcal{X}}$ as the \emph{partial extinction probability vector}. 
While global extinction implies partial extinction, there may be a positive chance that every type eventually disappears from the population while the total population size grows without bound; it is then possible that $\bm{q}<\bm{\tilde q}$ (see \cite[Section 5.1]{haut12} for an example).

At least partly due to these challenges, the set $S$ is yet to be fully characterised in the infinite-type setting.
There is, however, a number of papers that make progress toward this goal:
Moyal \cite{moyal62} gives general conditions for $S$ to contain at most a single solution $\vc s$ such that $\sup_{i\in \mathcal{X}} s_i <1$; Spataru \cite{Spa89} gives a stronger results by stating that $S$ contains at most two elements, $\bm{q}$ and $\bm{1}$; however, Bertacchi and Zucca
\cite{Zuc14,Zuc15} prove the inaccuracy of the latter by providing an irreducible example where $S$ contains uncountably many elements such that $\sup_{i\in \mathcal{X}} s_i=1$.
Both $\vc q$ and $\bm{\tilde{q}}$ are elements of the set $S$. It is well known that $\vc q$ is the minimal element, but as yet, there has been no attempt to identify the precise location of $\bm{\tilde q}$.
We observe that due to the existence of irreducible MGWBPs with $\bm{q}<\bm{\tilde q}<\bm{1}$, the partial extinction probability vector $\bm{\tilde q}$ may be neither the minimal element of $S$, which is $\bm{q}$, nor the maximal element of $S$, which is $\bm{1}$.

Extending the extinction criterion established in the finite-type case to the infinite-type setting has also proven difficult.
To resolve the problem in the infinite-type setting we should give both a partial and a global extinction criterion.
A number of authors have progressed in this direction
 \cite{Zuc14,Bra16,Har63,haut12,moyal62,Spa89,Zuc11}.
In the infinite-type case, the analogue of the Perron-Frobenius eigenvalue is the \emph{convergence norm} $\nu(M)$ of $M$ defined in \eqref{conv_norm}, which gives a partial extinction criterion: $\tilde{\bm{q}} = \bm{1}$ if and only if $\nu(M) \leq 1$, see \cite[Theorem 4.1]{Zuc11}.
However, when partial extinction is almost sure we are still lacking {general} necessary and sufficient conditions for $\vc q = \vc 1$. It turns out that 
there can be no global extinction criterion based solely upon $M$,
as highlighted through \cite[Example 4.4]{Zuc11}, but as pointed out by the author, other moment conditions have not been clearly identified.
In addition, when $\tilde{\bm{q}} < \bm{1}$, {following the terminology in \cite{Zuc14}}, the process can exhibit \textit{strong local survival} $\vc q=\tilde{\vc q}<\vc 1$, or \textit{non-strong local survival} $\vc q<\tilde{\vc q}<\vc 1$. It is again challenging to derive a general criterion separating the two cases.

The main contribution of this paper is to use a {unified} probabilistic approach to characterise the set $S$ and to derive a global extinction criterion applicable when $\tilde{\vc q}=\vc 1$
  for a class of branching processes with countably infinitely many types called \textit{lower Hessenberg branching processes} (LHBPs). In these processes, which have the typeset $\mathcal{X}= \{ 0,1,2, \dots \}$, the primary constraint is that \emph{type-$i$ individuals can produce offspring of type no larger than $i+1$}; as a consequence their (infinite) mean progeny matrices have a lower Hessenberg form. 
The probabilistic approach we employ relies on a single pathwise argument: we reduce the study of the LHBP to that of a much simpler Galton-Watson process in a varying environment (GWPVE), embedded in the LHBP. 
GWPVEs are single-type Galton-Watson processes whose offspring distributions vary deterministically with the generation.
In our context, the embedded GWPVE is \textit{explosive}, in the sense that individuals may have an infinite number of offspring. In particular, we show the equivalence between global extinction of the LHBP  and extinction of the embedded GWPVE, and between partial extinction of the LHBP and the event that all generations of the embedded GWPVE are finite. 
 Based on this relationship, we obtain several results for LHBPs:
\begin{itemize} 
  \item[\emph{(i)}] We prove that there is a continuum of fixed points solutions $\vc s \in S$, whose componentwise minimum and maximum are the global and partial extinction probability vectors $\bm{q}$ and $\bm{\tilde q}$, respectively (Theorem \ref{FPchar}).
\item[\emph{(ii)}] We establish a connection between the growth rates of the embedded GWPVE and the convergence rate of $s_i$ to $1$ as $i\to\infty$ for any $\vc s\in S\setminus\{\vc q,\vc 1\}$; this yields a physical interpretation for the fixed points lying in between $\vc q$ and $\tilde{\vc q}$  (Theorem \ref{RateDecay}).
\item[\emph{(iii)}] In the non-trivial case where $\tilde{\vc q}=\vc 1$, we provide a necessary and sufficient condition for global extinction which holds under some second moment conditions (Theorem \ref{KerstingExtCrit}). This is the first extinction criterion for irreducible processes that also applies to cases exhibiting non exponential growth. We illustrate the broad applicability of the criterion through some examples.
\item[\emph{(iv)}] Finally, under additional assumptions, we build on the global extinction criterion to derive necessary and sufficient conditions for strong local survival (Theorem \ref{Case4}).
\end{itemize}
While there is a {vast} literature on GWPVEs, the explosive case, 
 which {has already been studied} for standard Galton-Watson processes \cite{Lin15,Sag16}, is yet to be considered in the context of varying environment.
In order to prove our main theorems, we both apply known results on GWPVEs and develop new ones.
On the way to studying properties of the embedded GWPVE, we also derive a new partial extinction criterion for LHBPs which is computationally more efficient than other existing criteria.

The paper is organised as follows. In Section~2 we define LHBPs and introduce the tools we use to study them. In Section~3 we construct the embedded GWPVE and derive relationships between it and its corresponding LHBP. In Section~4 we develop \emph{(i)} and \emph{(ii)}. In Section 5 we deal with \emph{(iii)}. In Section 6 we illustrate the results of Section 5 through two examples. In Section 7 we address \emph{(iv)}. Finally, in Section 8 we discuss possible extensions of our results. 

In this paper, we let $\vc 1$ and $\vc 0$ denote the column vectors of $1$'s and $0$'s, respectively, and we let $\vc e_i$ represent the vector with all entries equal to zero, except entry $i$ which is equal to 1, the size of these vectors being defined by the context.  For any vectors $\vc x$ and $\vc y$, we write $\vc x\leq \vc y$ if $x_i\leq y_i$ for all $i$, and $\vc x< \vc y$ if $\vc x\leq  \vc y$ with $x_i<y_i$ for at least one entry $i$.

\section{Preliminaries}
Consider a MGWBP with the type set $\mathcal{X}=\mbN_0:=\{0,1,2, \dots \}$. We assume that the process initially contains a single individual whose type is denoted by $\varphi_0$. The process then evolves according to the following rules:
\begin{itemize}
\item[\it{(i)}] each individual lives for a single generation, and
\item[ \it{(ii)}] at death individuals of type $i$ give birth to $\bm{r}=(r_{  \ell})_{{  \ell}\in \{ 0,1, \dots, i+1\}}$ offspring, that is, $r_{0}$ individuals of type $0$, $r_{1}$ individuals of type $1$, \dots, and $r_{i+1}$ individuals of type $i+1$, where the vector $\bm{r}$ is chosen independently {of that  of all} other individuals according to a probability distribution, $p_i( \cdot)$, specific to the parental type $i \in \mathcal{X}$.
\end{itemize}
We refer to this as a \emph{lower Hessenberg branching process} (LHBP).

We construct the LHBP on the Ulam-Harris space \cite[Ch. VI]{Har63}, labelled $( \Omega, \mathcal{F}, \mbP)$, as follows. Let $\mathcal{J} = \bigcup_{n \geq 0} \mathcal{J}_n$ where $\mathcal{J}_n$ describes the virtual $n$-th generation. That is, $\mathcal{J}_0= \mathcal{X}$, where $\varphi_0 \in \mathcal{J}_0$ specifies the type of the root, and for $n\geq 1$, $\mathcal{J}_n={  \mathcal{X} \times} ( \mbN \times \mathcal{X} \times \mbN)^n$, where $(\varphi_0;i_1, j_1, y_1; \dots; i_n, j_n, y_n)$ denotes the $i_n$-th child of type $j_n$ born to $(\varphi_0; i_1, j_1, y_1; \dots; i_{n-1}, j_{n-1}, y_{n-1})$ and $y_n$ denotes the individual's unique identification number. Each virtual individual $I \in \mathcal{J}$ is assigned a random offspring vector $\bm{N}(I) =({N}_{\ell}(I))_{{  \ell} \in \mathcal{X}}$  that takes values in $R_j:=\{ \bm{r} \in (\mbN_0)^{\mathcal{X}}: r_\ell =0 \,\, \forall \, \, \ell > j+1\}$ when $I$ is of type $j$ and has distribution
 ${p}_j(\cdot)$, independently of  all other individuals. The random set of individuals who appear in the population, $X = \bigcup_{n \geq 0} X_n$, is then defined recursively from the values of $\bm{N}(I)$ as follows
\begin{equation}\label{Xd}
X_0 = \{ \varphi_0 \}, \quad X_n = \{ x = ( \tilde{x}; i_n, j_n, n) \in \mathcal{J}_n : \tilde{x} \in X_{n-1}, \, i_n \leq N_{j_n}(\tilde{x}) \}.
\end{equation}
The population in generation $n$ is described by the vector $\bm{Z}_n$ with entries 
$$
Z_{n,j} = \sum_{I \in \mathcal{J}_n }  \mathds{1} ( I \in X_n , \, \, j_n = j ), \quad j \in \mathcal{X}.
$$
We will often refer to branching processes by their sequence of population vectors $\{ \bm{Z}_n \}_{n \geq 0}$.

We define the \emph{progeny generating vector} $\bm{G}{ (\cdot)} : [0,1]^{\mathcal{X}} \to [0,1]^{\mathcal{X}}$, where
\begin{equation}\label{Gs}
G_i ( \bm{s} ) ={G_i (s_0,s_1,\ldots,s_{i+1} )} = \sum_{ \bm{r} \in R_i} p_{i}(\bm{r}) \bm{s}^{\bm{r}} = \sum_{ \bm{r} \in R_i } p_{i}(\bm{r}) \prod_{k=0}^{ i+1} s_k^{r_k},
\end{equation}
and the \emph{mean progeny matrix} $M=(M_{i,j})_{i,j\in\mathcal{X}}$, where
\begin{equation}\label{MPM}
M_{i,j} = \left. \left(\frac{\partial G_i ( \bm{s} )}{\partial s_j }\right)  \right|_{\bm{s}= \bm{1}}
\end{equation}
 is the expected number of type-$j$ children born to a parent of type $i$.
By assumption, $M$ is an infinite lower Hessenberg matrix.
To avoid trivialities we assume that $M_{i,i+1} >0$ for all $i \in \mathcal{X}$.
To $M$, we associate a weighted directed graph, referred to as the \emph{mean progeny representation graph}. 
This graph has vertex set $\mathcal{X}$ and contains an edge from $i$ to $j$ of weight $M_{i,j}$ if and only if $M_{i,j}>0$. 
The branching process is said to be \emph{irreducible} if there is a path between any two  vertices in the mean progeny representation graph on $\mathcal{X}$.
We define the \emph{convergence norm} of $M$,
\begin{equation}\label{conv_norm}
\nu(M)=\limsup_n \sqrt[n]{(M^n)_{ij}},
\end{equation}
which, when the process is irreducible, is independent of $i$ and $j$. 

The global and partial extinction probability vectors $\bm{q}$ and $\bm{\tilde q}$, defined in \eqref{global_def} and \eqref{partial_def}, are both solutions to the fixed point equation $\bm{s}=\bm{G}(\bm{s})$, and are thus elements of the set $S$ defined in \eqref{S}.
This can be seen by conditioning on the children of the initial individual and then observing that the process becomes partially (globally) extinct if and only if the daughter processes of these children become partially (globally) extinct.
Moreover, following the standard arguments, we can prove that $\bm{q}$ is the componentwise minimal element of $S$ (see \cite[Theorem 3.1]{moyal62}).
By the lower Hessenberg assumption, $\bm{s} = \bm{G}(\bm{s})$ can be written as $s_i=G_i(s_0,\ldots,s_i,s_{i+1})$ for all $i\geq 0$. Thus, by the monotonicity of $G_i(\cdot)$, each entry $s_i$ of any $\vc s\in S$ is uniquely determined by $s_0$.
It is then natural to consider the one-dimensional projection sets of $S$,
$$
S_i = \{ x \in [0,1] : \exists \, \bm{s} \in S, \text{ such that } s_i=x \},\quad {i\in\mathcal{X}}.
$$

\tikzset{
  treenode/.style = {align=center, inner sep=0pt, text centered,
    font=\sffamily},
  arn_n/.style = {treenode, circle, black, font=\sffamily\bfseries, draw=black,
    fill=white, text width=1.5em},
  arn_r/.style = {treenode, circle, red, draw=red, 
    text width=0em, very thick},
  arn_x/.style = {treenode, circle, white, font=\sffamily\bfseries, draw=black,
    fill=black, text width=1.5em}
}

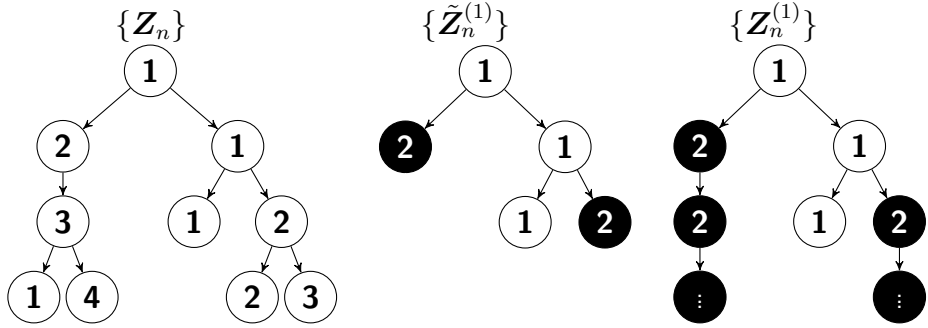
\begin{figure}
\begin{center}
$\{\bm{Z}_n\}$ \hspace{2.8cm} $\{\tilde{\bm{Z}}^{(1)}_n\}$ \hspace{2.6cm} $\{{\bm{Z}}^{(1)}_n\}$ \\
\begin{tikzpicture}[->,>=stealth',level/.style={sibling distance = 2.3cm/#1,
  level distance = 1cm}] 
\node [arn_n] {1}
    child{ node [arn_n] {2} 
            child{ node [arn_n] {3}  
            	child{ node [arn_n]{1}} 
		child{ node [arn_n]{4}}          								
            }                          
    }
    child{ node [arn_n] {1}
            child{ node [arn_n] {1}} 													            
            child{ node [arn_n] {2}
							child{ node [arn_n] {2}}
							child{ node [arn_n] {3}}
            }
		}
;  
\end{tikzpicture}
\quad
\begin{tikzpicture}[->,>=stealth',level/.style={sibling distance = 2.1cm/#1,
  level distance = 1cm}] 
\node [arn_n] {1}
    child{ node [arn_x] {2}
    	child{ node [white] {2}
		[white] child{ node [white]{\scalebox{1.5}{1}}}                          
    }
    }
    child{ node [arn_n] {1}
    	child{ node [arn_n] {1}} 
	child{ node [arn_x] {2}		
	} 
		}
;  
\end{tikzpicture}
\quad
\begin{tikzpicture}[->,>=stealth',level/.style={sibling distance = 2.1cm/#1,
  level distance = 1cm}] 
\node [arn_n] {1}
    child{ node [arn_x] {2} 
           child{ node [arn_x]{2}          
		child{ node [arn_x]{\scalebox{.55}{$\vdots$}}} 		    
           }                          
    }
    child{ node [arn_n] {1}
    	    child{ node [arn_n]{1}} 
             child{ node [arn_x]{2}          
		child{ node [arn_x]{\scalebox{.55}{$\vdots$}}} 		    
           } 
}
;  
\end{tikzpicture}
\end{center}
\caption{The processes $\{\bm{Z}_n \}$, $\{\bm{\tilde{Z}}^{(1)}_n\}$ and $\{\bm{{Z}}^{(1)}_n\}$ for a specific $\omega \in \Omega$.}
\label{RF2}
\end{figure}

We define two sequences of finite-type branching processes on $(\Omega, \mathcal{F}, \mbP)$. 
The  first, $\{ \bm{\tilde Z}^{(k)}_n \}_{n \geq 0, k \geq -1}$, is such that the random offspring vector of any virtual individual $I \in \mathcal{J}$ is given by
\[
\bm{\tilde N}^{(k)}(\omega, I) =
\begin{cases}
\bm{N}(\omega, I),\quad & t(I) \leq k \\
\bm{0}, \quad & t(I) > k, 
\end{cases}
\]
for any $\omega \in \Omega$, where $t(I)$ is the type of virtual individual $I$.
For any $k \geq -1$, outcomes of $\{ \bm{\tilde Z}^{(k)}_n \}$ are thus constructed by taking the corresponding outcome of $\{ \bm{Z}_n \}$ and removing the descendants of all individuals of type $i >k$. These types are said to be \emph{sterile}.
The second, $\{ \bm{Z}^{(k)}_n \}_{n \geq 0, k \geq -1}$, is such that the random offspring vector of any virtual individual $I \in \mathcal{J}$ is given by
\[
\bm{N}^{(k)}(\omega, I) =
\begin{cases}
\bm{N}(\omega, I),\quad & t(I) \leq k \\
\bm{e}_{t(I)}, \quad & t(I) > k,
\end{cases}
\]
for any $\omega \in \Omega$.
For any $k \geq -1$, outcomes of $\{ \bm{Z}^{(k)}_n \}$ are thus constructed by taking the corresponding outcome of $\{ \bm{Z}_n \}$ and replacing the descendants of all individuals of type $i >k$ with an infinite string of type-$i$ descendants. These types are said to be \emph{immortal}.
An illustration of $\{\bm{Z}_n \}$, $\{ \bm{\tilde Z}^{(1)}_n \}$ and $\{ \bm{Z}^{(1)}_n \}$ for a specific $\omega \in \Omega$ is given in Figure \ref{RF2}.
By construction, for all $\omega\in\Omega$, \emph{(i)} for each fixed value of $k$, if $\varphi_0 \leq k+1$ then the sterile and immortal individuals are necessarily of type $k+1$, \emph{(ii)} 
\begin{equation}\label{mod equiv} 
Z_{n,\ell}^{(k)}(\omega) = \tilde{Z}_{n,\ell}^{(k)}(\omega)\qquad\textrm{for all $n \geq 0$ and $0 \leq \ell \leq k$,}
\end{equation}
and \emph{(iii)}
\begin{equation}\label{mod equiv2} 
\lim_{n \to \infty} {Z}_{n,k+1}^{(k)}(\omega) = \sum_{n=0}^\infty \tilde{Z}^{(k)}_{n,k+1}(\omega).
\end{equation}
We denote the progeny generating vector of $\{ \bm{Z}^{(k)}_n \}$ by $\bm{G}^{(k)}( \bm{s})$, which has entries
\begin{equation}\label{modGF}
{G}_i^{(k)}( \bm{s}) = 
\begin{cases}
G_{i}( \bm{s} ), \quad &  0 \leq i \leq k \\
s_{k+1}, \quad & i =k+1.
\end{cases}
\end{equation}
By Equation \eqref{mod equiv}, the global extinction probability vectors of $\{ \bm{Z}_n^{(k)} \}$ and $\{ \bm{\tilde{Z}}^{(k)}_n \}$, denoted by $\bm{\tilde{q}}^{(k)}$ and $\bm{q}^{(k)}$, are given by
\[
\bm{\tilde{q}}^{(k)} = \lim_{n \to \infty} \bm{G}^{(k,n)} ( 0, \dots, 0, 1) \quad \mbox{ and } \bm{{q}}^{(k)} = \lim_{n \to \infty} \bm{G}^{(k,n)} ( 0, \dots, 0, 0),
\]
where $\bm{G}^{(k,n)}(\cdot)$ is the $n$-fold composition of $\bm{G}^{(k)}(\cdot)$. 
As demonstrated in \cite[Lemma 3.1]{haut12}, the sequence $\{\bm{q}^{(k)}\}_{k\geq -1}$ increases componentwise to $\bm{q}$. In addition, if $\{ \bm{Z}_n \}$ is irreducible and non-singular (that is, there exists $i \in \mathcal{X}$ such that $\sum _{\bm{v}: |\bm{v}|=1} p_{i}(\bm{v}) <1$), then the sequence $\{\bm{\tilde{q}}^{(k)}\}_{k\geq -1}$ decreases componentwise to $\bm{\tilde{q}}$ (see Theorem~\ref{PartConv} in Appendix A). Note that \cite[Lemma 3.2]{haut12} is an inaccurate version of the latter result. 
Unless stated otherwise,  we assume that $\{ \bm{Z}_n \}$ is non-singular and irreducible.

\section{An embedded GWPVE with explosions}\label{sec_GWPVE}

\tikzset{
  treenode/.style = {align=center, inner sep=0pt, text centered,
    font=\sffamily},
  arn_n/.style = {treenode, circle, black, font=\sffamily\bfseries, draw=black,
    fill=white, text width=1.5em},
  arn_r/.style = {treenode, circle, red, draw=red, 
    text width=1.5em, very thick},
  arn_x/.style = {treenode, circle, white, font=\sffamily\bfseries, draw=black,
    fill=black, text width=1.5em}
}

\begin{figure}
\centering
\hspace{.5cm}$\{ \bm{Z}_n \}$ \hspace{5.3cm} $\{Y_k\}$ \\
\begin{tikzpicture}[->,>=stealth',level/.style={sibling distance = 2.4cm/#1,
  level distance = 1cm}] 
\node [arn_x] {0}
    child{ node [arn_n] {0} 
    	child{ node [arn_n] {0}}
	child{ node [arn_x] {1}
		child{ node [arn_n] {1}
			child{ node [arn_x] {2}
				child{ node [arn_x] {3}}
			}			
		}
		child{ node [arn_x] {2}
			child{ node [arn_n] {1}}
		}
	}                    
    }
    child{ node [arn_x] {1}
            child{ node [arn_n] {0}
		child{ node [arn_n]{1}} 						
            }
            child{ node [arn_x] {2}
		child{ node [arn_n] {2}}
		child{ node [arn_x] {3}
			child{ node [arn_n] {3}
				child{ node [arn_x] {4}}
			}
		}
            }     
    }
;  
\end{tikzpicture} \hspace*{2cm}
\begin{tikzpicture}[->,>=stealth',level/.style={sibling distance = 2.1cm/#1,
  level distance = 1cm}] 
\node [arn_x] {0}
		child{ node [arn_x]{1}
			child{ node [arn_x]{2}
				child{ node [arn_x]{3}
					child{ node [white] {1}
						[white] child{ node [white]{\scalebox{1.5}{0}}}                          
    }
				}
			}
			child{ node [arn_x]{2}}
		}
		child{ node [arn_x]{1}
			child{ node [arn_x]{2}
				child{ node [arn_x]{3}
					child{ node [arn_x]{4}}
				}
			}
		}

;    
\end{tikzpicture}
\caption{An outcome of $\{ \bm{Z}_n \}$ and $\{Y_k\}$ for a specific $\omega \in \Omega$. The {highlighted} type-$k$ individuals 
represent the sterile individuals in the corresponding realisation of $\{ \bm{\tilde{Z}}^{(k-1)}_n \}$. 
}
\label{embedvar}
\end{figure}
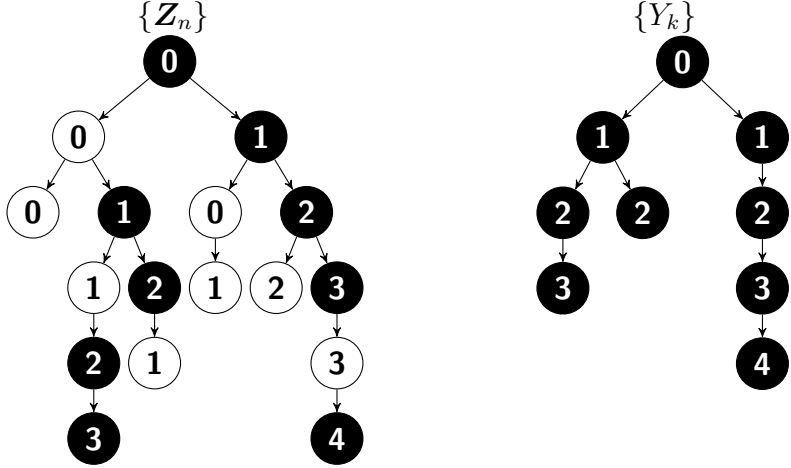

We construct the embedded GWPVE $\{ Y_k \}$ on $( \Omega, \mathcal{F}, \mbP)$ from the paths of ${\{ \bm{{Z}}_n \}}$ by selecting all individuals whose type is strictly larger than that of all their ancestors, and connecting each selected individual to their nearest (in generation) selected ancestor (see Figure \ref{embedvar}). 
More formally, we define a function $f(\cdot) : \mathcal{J} \to \mathcal{J}$ that takes a line of descent $(\varphi_0; i_1, j_1,  y_1;  \dots; i_n, j_n,   y_n)$ and deletes each  triple $(i_k, j_k,  y_k)$ whose type is not strictly larger than all its ancestors.
For each $\omega \in \Omega$ the family tree of $\{{Y}_k \}$ is then given by $f(X(\omega))$, where $X(\omega)$ is defined in \eqref{Xd}.
Variants of $\{ Y_k \}$ (which do not permit explosion) can be found in \cite{Bra16} and \cite{Gan10}.

We take the convention that $\{ Y_k\}$ starts at the generation number corresponding to the initial type $\varphi_0$ in $\{ \bm{Z}_n \}$. By construction, for any $\omega\in\Omega$ we then have
\begin{equation}\label{EmbeddedGWPVE}
Y_k(\omega)= \sum^\infty_{n=0} \tilde{Z}^{({ k-1})}_{n,{ k}}(\omega) = \lim_{n \to \infty} {Z}^{({ k-1})}_{n,{ k}} ( \omega),
\end{equation}
that is, the $k$th generation of $\{Y_k\}$ is made up every sterile (type-$k$) individual produced over the lifetime of $\{ \bm{\tilde{Z}}_n^{({k-1})} \}_{n\geq 0}$.
By the lower Hessenberg assumption, each sterile type-$k$ individual that appears in $\{ \bm{\tilde Z}^{(k-1)}_n \}_{n\geq 0}$ is a descendant of a sterile type-$(k-1)$ individual that appears in $\{ \bm{\tilde Z}^{(k-2)}_n \}_{n\geq 0}$.
Thus, because the daughter processes of these type-($k-1$) individuals in $\{ \bm{\tilde Z}_n^{(k-1)} \}$ are i.i.d., $Y_k$ satisfies the branching process equation
\begin{equation}\label{P2BPeq}
Y_k \stackrel{d}{=} \sum^{Y_{k-1}}_{i=1} \xi_{k,i},
\end{equation}
where $\{ \xi_{k,i} \}_{i \geq 1}$ is a sequence of i.i.d. random variables such that, $\xi_{k,i} \stackrel{d}{=} \sum^\infty_{n=0} \tilde Z_{n,k}^{(k-1)}$  conditional on $\varphi_0=k-1$.
This means $\{ Y_k \}$ is indeed a GWPVE;
it is however not a classical one because $\{ Y_k \}$ may have a positive chance of \emph{explosion}, that is, individuals in $\{ Y_k \}$ may give birth to an infinite number of offspring with positive probability. 
The next lemma states that $\{ Y_k \}$ explodes by generation $k$ if and only if $\{ \bm{\tilde{Z}}_n^{({k-1})} \}_{n\geq 0}$ survives globally, and $\{ Y_k \}$ becomes extinct by generation $k$ if and only if $\{ \bm{{Z}}_n^{({k-1})} \}_{n\geq 0}$ becomes globally extinct.

\begin{lemma}\label{P2L1}
For any $k \geq \varphi_0$,
\begin{equation}\label{fa}
\{ \omega \in \Omega : Y_k(\omega) < \infty \} \stackrel{a.s.}{=} \{ \omega \in \Omega : \lim_{n \to \infty}  \bm{\tilde Z}_n^{(k-1)}(\omega)=\bm{0} \},\end{equation}
and
\begin{equation}\label{sa}
\{ \omega \in \Omega :  Y_k(\omega) = 0 \} \stackrel{a.s.}{=}  \{ \omega \in \Omega : \lim_{n \to \infty} \bm{Z}_n^{(k-1)}(\omega)=\bm{0} \}.
\end{equation}
\end{lemma}

\noindent\textbf{Proof:}
To prove \eqref{fa}, first suppose that $\omega \in \{\lim_n \bm{\tilde Z}^{(k-1)}_n = \bm{0} \}$. Then there exists a generation $N<\infty$ such that $\bm{\tilde Z}^{(k-1)}_n(\omega)=\bm{0}$ for all $n \geq N$. This means $\sum_{n=0}^\infty {\tilde{Z}}^{(k-1)}_{n,k}(\omega) =\sum_{n=0}^N  {\tilde{Z}}^{(k-1)}_{n,k}(\omega)  < \infty$, which implies $\omega \in \{Y_k < \infty \}$.
It then remains to prove 
$\mbP( Y_k < \infty ,\, \, \liminf_n |\bm{\tilde Z}^{(k-1)}_n| > {0} )=0.$
Because $M_{i,i+1}>0$ for any $i \leq k-1$, we have 
$
\mbP_i \left( \sum_{n=0}^{k} \tilde{Z}^{(k-1)}_{n,k} =0\right) \leq 1-\varepsilon_i,
$
for some $\varepsilon_i>0$.
Thus, for any  {$\bm{\tilde{z}}^{(k-1)}_0 \in (\mbN_0)^{k+1}$},  we have
\begin{equation*}
\mbP \left( \sum_{n=0}^{k} \tilde{Z}^{(k-1)}_{n,k} =0, \, |\bm{\tilde{Z}}^{(k-1)}_{k}|>0 \, \Big| \, \bm{\tilde{Z}}^{(k-1)}_0=\bm{\tilde{z}}^{(k-1)}_0\right) \leq 1-\varepsilon.
\end{equation*}where $\varepsilon := \min_{0\leq i \leq k} \{ \varepsilon_i\} >0$.
By the Markov property, we then have
\begin{align*}
\mbP( &Y_k< \infty, \, \liminf_{n \to \infty} |\bm{\tilde{Z}}_n^{(k-1)} |> {0} ) \\
 &\leq \mbP\left(\exists N{\geq 0} : \, \sum_{{ m}=N}^{\infty} \tilde{Z}_{m,k}^{(k-1)}=0, \, | \bm{\tilde{Z}}_n^{(k-1)}| > 0 \, \forall n\right)\\
 &= \mbP\left( \exists N\geq 0 : \, { \bigcap_{\ell=0}^\infty \left( \sum_{m=N+\ell k}^{N+(\ell+1)k} \tilde{Z}_{m,k}^{(k-1)}=0, \, | \bm{\tilde{Z}}^{(k-1)}_{ N+(\ell+1)k} |>0\right)}\right)  \\
&{\leq \sum^\infty_{N=0}\prod_{\ell \geq 0} (1-\varepsilon)=0},
\end{align*} 
leading to \eqref{fa}. The same arguments lead to \eqref{sa}.
\qed \medskip

Since
$$\{ \omega \in \Omega : \liminf_n \bm{\tilde Z}_n^{(k-1)}(\omega) > \bm{0} \} \subseteq \{ \omega \in \Omega: \liminf_n \bm{\tilde Z}_n^{(k)}(\omega) > \bm{0} \},$$ and $$\{ \omega \in \Omega : \lim_n \bm{ Z}_n^{(k-1)}(\omega) = \bm{0} \} \subseteq \{ \omega \in \Omega: \lim_n \bm{ Z}_n^{(k)}(\omega) = \bm{0}\},$$ 
the process $\{ Y_k \}$ has two absorbing states, $0$ and $\infty$.
The next corollary formalises the equivalence between the following
events:
$$\begin{array}{ccc} \{ \bm{Z}_n \} \mbox{ experiences:} & &\{ Y_k \} \mbox{ reaches:}\\
\mbox{both partial and global extinction} & \equiv & \mbox{the absorbing state $0$}\\\mbox{neither partial nor global extinction} & \equiv & \mbox{the absorbing state $\infty$}\\\mbox{partial extinction but not global extinction} & \equiv & \mbox{neither $0$ nor $\infty$.}

\end{array}$$

\begin{corollary}\label{P2C1}
The global extinction event $\mathcal{E}_g \stackrel{a.s.}{=} \{ \omega \in \Omega : \lim_{k\to\infty} Y_k(\omega) =0 \}$, and the partial extinction event $\mathcal{E}_p \stackrel{a.s.}{=} \{ \omega \in \Omega : Y_k(\omega) < \infty, \, \forall k \geq \varphi_0\}$.
\end{corollary}
\noindent\textbf{Proof:}
The result follows from Lemma \ref{P2L1} and the arguments in the proofs of \cite[Lemmas 3.1]{haut12} and Theorem~\ref{PartConv} respectively. 
\qed \medskip

By Corollary \ref{P2C1} we can express  any question about the extinction probability vectors $\bm{q}$ and $\bm{\tilde q}$ in terms of the process $\{ Y_k \}$. 
In the sequel we use the shorthand notation $\mbP_i(\cdot)$ for $\mbP(\cdot | Y_i=1)$ and $\mbE_i(\cdot)$ for $\mbE(\cdot | Y_i=1)$.

\begin{corollary}\label{Link}
For any $k \geq 0$ and $0\leq i\leq k$,
$$
q_i^{(k)}=\mbP_i\left( Y_{{k+1}}=0 \right) \quad \mbox{ and } \quad \tilde{q}_i^{(k)}= \mbP_i\left( Y_{{k+1}}< \infty \right), 
$$
and for any $i \geq 0$,
$$
q_i=\mbP_i(\lim_{k \to \infty} Y_k =0) \quad \mbox{ and } \quad \tilde{q}_i = \mbP_i \left( \forall \, k \geq {i}, \, Y_k < \infty \right).
$$
\end{corollary}
\noindent\textbf{Proof:}
The results are immediate consequences of Lemma \ref{P2L1} and Corollary~\ref{P2C1}.
\qed \medskip 

To take advantage of Corollary \ref{Link} we require the progeny generating function of each generation of the embedded GWPVE.
For $k \geq 0$, we let
 \begin{equation}\label{offspring} g_k(s):=\mbE_k(s^{Y_{k+1}}\,\mathds{1}\{Y_{k+1}<\infty\})=\sum_{\ell\geq 0} \mbP \left( \sum^\infty_{n=1} \tilde{Z}^{(k)}_{n,k+1} = \ell \,\Big|\,\varphi_0=k\right) s^\ell ,\end{equation} where $s \in [0,1] $.
Due to the possibility of explosion we may have $1>g_k(1)=\mbP_k(Y_{k+1}<\infty)=\tilde{q}_k^{(k)}$.
By \eqref{P2BPeq}, the generating function of $Y_{k+1}$, conditional on $Y_i=1$ for $i \leq k$, is given by
\[
g_{i \to {k}}(s):=g_i \circ g_{i+1} \circ \dots \circ g_{k}(s), \quad s \in [0,1].
\]
Consequently, by Corollary \ref{Link}, we have $q^{(k)}_i= g_{i \to {k}} (0)$, $ \tilde{q}^{(k)}_i= g_{i \to {k}} (1)$,
$q_i = \lim_{k \to \infty} g_{i \to k} (0)$, and $\tilde{q}_i= \lim_{k \to \infty} g_{i \to k} (1).$

The next two lemmas provide respectively an explicit and an implicit relation between the sequence of progeny generating functions $\{g_k(\cdot)\}$ and the progeny generating vector $\vc G(\cdot)$.  The first requires the following technical assumption:

\begin{assumption}\label{Nonsink}
For all $k\geq{0}$, 
\begin{equation}\label{truncdicS}
\mbP_k \left( \lim_{n \to \infty} \sum_{i =0}^k Z_{n,i}^{(k)} \to 0\right)+\mbP_k \left(\lim_{n \to \infty}  \sum_{i =0}^k Z_{n,i}^{(k)} \to \infty \right)=1.
\end{equation}
\end{assumption}

\begin{lemma}\label{EGF} 
If Assumption \ref{Nonsink} holds,
then for all $k \geq {\varphi_0}$, the progeny generating function of $\{ Y_k \}$ at generation $k$ is given by 
\[
g_k ( s ) = \lim_{n \to \infty} G_k^{(k,n)} ( s_{0}, s_{1}, \dots, s_k, s) , \quad s \in [0,1],
\]
where $(s_{0}, s_{1}, \dots, s_k ) \in [0,1)^{k+1}$.
\end{lemma}

\noindent\textbf{Proof:}
By \eqref{EmbeddedGWPVE} and \eqref{offspring}, 
\begin{align}\label{EGFEE1}
g_k(s) &= \mbE_k \left( s^{\lim_{n \to \infty} Z_{n, k+1}^{(k)} } \,\mathds{1}\left\{ \lim_{n \to \infty} Z_{n,k+1}^{(k)} < \infty \right\}  \right).
\end{align}
By Assumption \ref{Nonsink} and the fact that $(s_0, \dots, s_k) \in [0,1)^{k+1}$, 
\[
\mbP_k \left( \lim_{n \to \infty} \prod_{{ i=0}}^k s_i^{Z^{(k)}_{n,i}}=0 \right) + \mbP_k \left( \lim_{n \to \infty} \prod_{{ i=0}}^k s_i^{Z^{(k)}_{n,i}}=1\right)=1,
\]
that is, $\lim_{n \to \infty} \prod_{{ i=0}}^k s_i^{Z^{(k)}_{n,i}}$ is an indicator function. In addition, Lemma \ref{P2L1} implies
$
\{  \lim_{n \to \infty} Z_{n,k+1}^{(k)} < \infty\} \stackrel{a.s.}{=} \{  \lim_{n \to \infty} \bm{\tilde Z}_n^{(k)} = \bm{0} \} 
=  \{ \lim_{n \to \infty} \prod_{{ i=0}}^k s_i^{\tilde Z^{(k)}_{n,i}}=1\} = \{ \lim_{n \to \infty} \prod_{{ i=0}}^k s_i^{Z^{(k)}_{n,i}}=1\}.
$
Thus, \eqref{EGFEE1} can be rewritten as 
$$
g_k ( s ) = \mbE_k \left(\lim_{n \to \infty} s^{Z^{(k)}_{n,k+1}}  \prod_{{i=0}}^k s_i^{Z^{(k)}_{n,i}} \right) =\lim_{n \to \infty} G_k^{(k,n)} ( s_{0}, s_{1}, \dots, s_k, s), $$
where the last equality follows from the dominated convergence theorem. 
\qed \medskip

\begin{lemma}\label{Genlink}
For {any $k\geq 0$, the progeny generating function $g_{k}(\cdot)$ satisfies}
\begin{equation}\label{VEBPtoITBP}
g_{k}(s)= {G}_k \left( g_{0 \to k}(s), g_{1 \to k}(s), \dots,  g_k(s), s  \right).
\end{equation}
\end{lemma}

\noindent\textbf{Proof:}
By conditioning on the offspring of a type-$k$ individual in $\{\vc{\tilde{Z}}^{(k)}_{n}\}$,
\small{\begin{eqnarray*}\lefteqn{ g_k(s)}\\&=&\mbE\left[ s^{\sum_{n=1}^\infty  \tilde{Z}^{(k)}_{n,k+1}} \,  \mathds{1}\left\{\scriptstyle\sum_{n=1}^\infty  \tilde{Z}^{(k)}_{n,k+1} < \infty \right\}\,\Big|\, \varphi_0=k\right] \\&=&\sum_{\vc z\geq \vc 0} \mbE\left[ s^{\sum_{n=1}^\infty  \tilde{Z}^{(k)}_{n,k+1}} \,  \mathds{1}\left\{\scriptstyle\sum_{n=1}^\infty  \tilde{Z}^{(k)}_{n,k+1} < \infty \right\}\,\Big|\, \varphi_0=k, \vc{\tilde{Z}}^{(k)}_{1}=\vc z \right] \,\mbP[\vc{\tilde{Z}}^{(k)}_{1}=\vc z\,|\,\varphi_0=k]. 
\end{eqnarray*}} 
\normalsize Then, by the Markov property and the independence between the daughter processes of individuals from the same generation, 
\begin{eqnarray}
\lefteqn{  \mbE\left[ s^{\sum_{n=1}^\infty  \tilde{Z}^{(k)}_{n,k+1}} \,  \mathds{1}\left\{\scriptstyle\sum_{n=1}^\infty  \tilde{Z}^{(k)}_{n,k+1} < \infty \right\}\,\Big|\, \varphi_0=k, \vc{\tilde{Z}}^{(k)}_{1}=(z_0,\ldots,z_k,z_{k+1}) \right] } \nonumber \\
&=&s^{z_{k+1}} \prod_{i=0}^k\mbE\left[ s^{\sum_{n=1}^\infty  \tilde{Z}^{(k)}_{n,k+1}} \,  \mathds{1}\left\{\scriptstyle\sum_{n=1}^\infty  \tilde{Z}^{(k)}_{n,k+1} < \infty \right\}\,\Big|\, \varphi_0=i\right]^{z_i} \nonumber \\ &=&s^{z_{k+1}} \prod_{i=0}^k g_{i\rightarrow k}(s)^{z_i}, \label{SLFe}
\end{eqnarray} 
where \eqref{SLFe} follows from \eqref{P2BPeq}. 
This leads to
$$ 
g_k(s) = \sum_{\vc z\geq \vc 0}  \prod_{i=0}^k g_{i\rightarrow k}(s)^{z_i} \, s^{z_{k+1}}\,\mbP[\vc{\tilde{Z}}^{(k)}_{1}=\vc z\,|\,\varphi_0=k],
$$
which completes the proof. 
\qed \medskip

\section{Fixed points and extinction probabilities}\label{fpep}

We now characterise the set $S$ defined in \eqref{S}. The main results in this section rely on the relation between $S$ and the set
$$
S^{[e]}=\{ \bm{s}\in[0,1]^{\mathcal{X}} : s_{k} = {g}_k ( {s}_{k+1} ) \, \forall \, k \geq 0 \},
$$
which corresponds to the set of fixed points of the embedded GWPVE. Because each ${g}_k(\cdot)$ is a monotone increasing function, like $S$, the set $S^{[e]}$ is one-dimensional. 
In this section we assume that Assumption \ref{Nonsink} holds.
For any vector $\vc s\in S$, we write $\bm{\bar{s}}^{(k)}:=( s_0, s_1, \dots, s_{k} )$ for the restriction of $\vc s$ to its first $k+1$ entries.

The next lemma establishes a relationship between $S$ and $S^{[e]}$.

\begin{lemma}\label{FPequ} $S=S^{[e]} \cup \{\bm{1} \}.$
\end{lemma}

\noindent\textbf{Proof:}
Suppose $\bm{s} \in S$ and $\bm{s} \neq \bm{1}$.  
For any $k,n\geq 0$, $\bm{\bar{s}}^{(k+1)}$ satisfies $\bm{\bar{s}}^{(k+1)}= {\bm{G}}^{(k,n)}( \bm{\bar{s}}^{(k+1)})$.
Because $\{\vc Z_n\}$ is irreducible and $\bm{s} \neq \bm{1}$ we have $s_i <1$ for all $i \in \mathcal{X}$ (see \cite[Theorem 2]{Spa89}).
Thus, using Lemma~\ref{EGF},
$
g_k(s_{k+1}) = \lim_{n \to \infty} {G}_k^{(k,n)} (\bm{\bar{s}}^{(k+1)} ) = s_k
$
for all $k \geq 0$, leading to $\bm{s} \in S^{[e]}$.
Now suppose $\bm{s} \in S^{[e]}$. Then, by Lemma \ref{Genlink}, for all $k \geq 0$,
\begin{align*}
s_k &= g_k(s_{k+1}) ={G}_k \left( g_{0 \to k}(s_{k+1}), g_{1 \to k}(s_{k+1}), \dots,  g_k(s_{k+1}), s_{k+1}  \right) = {G}_k ( \bm{s} ),
\end{align*}therefore $\bm{s} \in S$.
\qed \medskip

We now characterise the one-dimensional projection sets $S_i$ and identify which elements of $S$ correspond to the global and partial extinction probability vectors.

\begin{theorem}\label{FPchar}If $S=\{\vc 1\}$ then $\bm{q}=\bm{\tilde{q}}=\bm{1}$; otherwise
$$\bm{q}={\min} \, S \quad \mbox{and} \quad \bm{\tilde{q}} = \sup S \backslash  \{\vc1\} .$$ 
In particular,
$S_i=[q_i, \tilde{q}_i]  \cup 1\, \text{ for all $i\geq 0$.}$
\end{theorem}

\noindent\textbf{Proof:}
We show that
\begin{equation}\label{infsup}\bm{q}=\min S^{[e]} \quad \mbox{and} \quad \bm{\tilde{q}} = \max{S^{[e]}},\end{equation}
and {for any $i\geq 0$,}
$ S^{[e]}_i=[q_i, \tilde{q}_i]$,
where $$
S_i^{[e]} = \{ x \in [0,1] : \exists \, \bm{s} \in S^{[e]}, \text{ such that } s_i=x \}.
$$
These results follow from the fact that $g_i(\cdot)$ and $g^{-1}_i(\cdot)$ are monotone increasing functions, and therefore so are $g_{i \to j}(\cdot)$ and $g^{-1}_{i \to j-1}(\cdot) := g^{-1}_{ j-1} \circ \dots \circ g^{-1}_i(\cdot)$ for $j>i$. Let $\bm{s} \in S^{[e]}$, then for all $0 \leq i <k$,
\[
{q_i^{ (k-1)} = } g_{i \to k-1}(0) \leq s_i=g_{i \to k-1}(s_k) \leq g_{i \to k-1}(1) =\tilde{q}^{({k-1})}_i.
\]
Taking the limit as $k \to \infty$ we obtain $q_i \leq s_i \leq \tilde{q}_i$ for all $i \geq 0$, which shows \eqref{infsup}. 
Now suppose $q_i \leq s_i \leq\tilde{q}_i$. 
For any $ j <i$, define $s_j:= g_{j \to i-1}(s_i)$; then
\[
q_j = g_{j \to i-1}(q_i) \leq {s_j} \leq g_{j \to i-1}(\tilde{q}_i)= \tilde{q}_j.
\]
Similarly, for any $j>i$, define $s_j:= g^{-1}_{i \to {j-1}}(s_i) $; then
\[
q_j = g^{-1}_{i \to {j-1}}(q_i) \leq{s_j }\leq g^{-1}_{i \to { j-1}}(\tilde{q}_i)= \tilde{q}_j.
\]
This shows that for any $i\geq 0$ and for any $s_i \in [q_i, \tilde{q}_i]$, it is possible to construct a vector $\vc s$ belonging to $S^{[e]}$. 
\qed \medskip

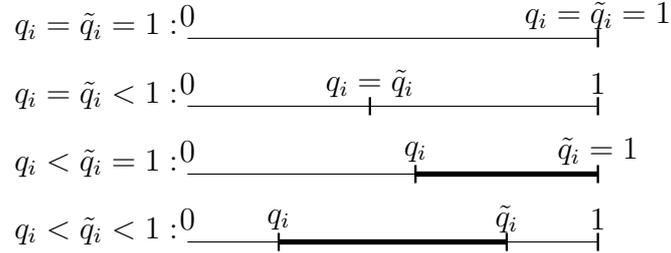
\begin{figure}
\centering
\begin{tikzpicture}[scale=0.6]

\node at (1,3.25){ ${q}_i={\tilde{q}}_i< {1}:$};
\node at (1,4.75){ ${q}_i={\tilde{q}}_i= {1}:$};
\node at (1,1.75){ ${q}_i<{\tilde{q}}_i= {1}:$};
\node at (1,0.25){ ${q}_i<{\tilde{q}}_i< {1}:$};

\draw (3,4.5) to (12,4.5);
\draw (3,3) to (12,3);
\draw (3,1.5) to (12,1.5);
\draw (3,0) to (12,0);

\draw[line width=0.3mm] (12,4.3) to (12,4.7);
\draw[line width=0.3mm] (12,2.8) to (12,3.2);
\draw[line width=0.3mm] (12,1.3) to (12,1.7);
\draw[line width=0.3mm] (12,-.2) to (12,.2);


\node at (3,5){$0$};
\node at (3,3.5){$0$};
\node at (3,2){$0$};
\node at (3,0.5){$0$};

\draw[line width=0.3mm] (7,2.8) to (7,3.2);

\draw[line width=0.7mm] (8,1.5) to (12,1.5);
\draw[line width=0.3mm] (8,1.3) to (8,1.7);

\draw[line width=0.7mm] (5,0) to (10,0);
\draw[line width=0.3mm] (5,-0.2) to (5,0.2);
\draw[line width=0.3mm] (10,-0.2) to (10,0.2);

\node at (7,3.5){$q_i=\tilde{q}_i$};
\node at (12,3.5){$1$};

\node at (12,5){$q_i=\tilde{q}_i=1$};

\node at (12,2){$\tilde{q}_i=1$};
\node at (8,2){${q}_i$};

\node at (12,0.5){$1$};
\node at (5,0.5){${q}_i$};
\node at (10,0.5){$\tilde{q}_i$};

\end{tikzpicture}
\caption{A visual representation of {possible} sets $S_i$ in the irreducible case.} 
\label{Fixed Vis} 
\end{figure}

Theorem \ref{FPchar} implies that $S$ contains one, two, or uncountably many elements. More specifically, it shows that $\bm{q}$ is the minimal element of $S$
 which is the beginning of a continuum of elements whose supremum is $\bm{\tilde{q}}$, as illustrated in Figure~\ref{Fixed Vis}.

\begin{remark}\label{rem1}In the reducible case there may be an additional countable number of fixed points $\vc s$ such that $\vc{\tilde q}\leq \vc s\leq \vc 1$. We refer to \cite[Section 4.4]{BraPhD} for the details.
\end{remark}

With the goal of giving a probabilistic interpretation to the intermediate fixed points $\vc s \in S$ such that $\bm{q} < \bm{s} < \bm{\tilde{q}}$,
we now derive properties of the infinite-dimensional set $S$.
We begin by deriving a sufficient condition for
\begin{equation}\label{MoyalEq}
\lim_{i \to \infty} s_i =1, \quad \text{ for all } \bm{s} \in S \backslash \{ \bm{q} \},
\end{equation}
that is, for $S$ to contain \emph{at most} a single element (corresponding to $\bm{q}$) whose entries do not converge to $1$. 
In a more general setting, sufficient conditions for \eqref{MoyalEq} can be found in Moyal \cite[Lemmas 3.3 and 3.4]{moyal62}, the most notable being `$\inf q_i > 0$'. The same author also conjectures a more general condition:
\[
\sup_i p_i^{(1)} <1, \quad \quad \text{ where } \quad p_i^{(1)} := \sum_{\bm{v}: |\bm{v}|=1} p_i(\bm{v}).
\]
In the case of LHBPs we now provide a stronger result.

\begin{theorem}\label{LimInfS}
If 
\begin{equation}\label{LimSCond}
\sum_{i=0}^\infty (1-p_i^{(1)}) = \infty
\end{equation}  
then \eqref{MoyalEq} holds.
\end{theorem}
The proof of Theorem \ref{LimInfS} uses the following lemma which we state separately because, for LHBPs, it generalises the conditions of \cite[Theorem~1]{Bra16}.

\begin{lemma}\label{LemYkdic}
If \eqref{LimSCond} holds then 
$\mbP(Y_k \to 0) + \mbP(Y_k \to \infty)=1$.
\end{lemma}
\noindent\textbf{Proof:}
Following Lindvall \cite{Lind74}, we have $\mbP(Y_k \to 0) + \mbP(Y_k \to \infty)=1$ if and only if 
\begin{equation}\label{EmGWPVEdicho}
\sum^\infty_{k=0} (1-g'_k(0)) = \infty.
\end{equation}
Suppose \eqref{LimSCond} holds without \eqref{EmGWPVEdicho}, that is, assume \eqref{LimSCond} and 
\begin{equation}\label{NEmGWPVEdicho}
\sum^\infty_{k=0} (1-g'_k(0)) < \infty.
\end{equation}
In this case, there can be only finitely many $k$ such that $g'_k(0)=0$. 
Thus, 
\eqref{NEmGWPVEdicho} holds if and only if there exists $\ell \geq 0$ such that 
\begin{equation}\label{NEmGWPVEdicho2}
\prod^\infty_{k=\ell} g'_k(0) \equiv \mbP_\ell ( Y_k =1, \, \forall k \geq \ell)>0.
\end{equation}
In addition, because $M_{i,i+1} >0$ for all $i\geq 0$, in every generation of the embedded process (including any for which $g'_k(0)=0$) individuals have a positive chance of giving birth to at least one offspring. In combination with \eqref{NEmGWPVEdicho2} this implies that there exists $c>0$ such that for any $l \geq 0$, 
$
\mbP_l ( Y_k \geq 1, \, \forall k \geq l) \geq c.$
Recall that each individual in $\{Y_k\}$ corresponds to an individual in $\{\bm{Z}_n\}$. If the corresponding individual in $\{ \bm{Z}_n \}$ has no offspring then neither does the individual in $\{Y_k \}$, whereas if the corresponding individual in $\{ \bm{Z}_n \}$ has two or more offspring then the individual in $\{Y_k\}$ must have at least two offspring with probability greater than or equal to $c^2$. Thus, for all $k\geq 0$,
$
1-g'_k(0) \geq c^2(1-p_k^{(1)}),$
which implies 
$
\sum_{k=0}^\infty (1-g'_k(0)) \geq c^2 \sum_{k=0}^\infty (1-p_k^{(1)}) =\infty,
$
contradicting \eqref{NEmGWPVEdicho}. Therefore, if \eqref{LimSCond} holds, we must have \eqref{EmGWPVEdicho}.
\qed \medskip

\noindent\textbf{Proof of Theorem \ref{LimInfS}:} 
By Lemma \ref{FPequ}, we may assume $\bm{s} \in S^{[e]}$.  Thus, for all
$k \geq 0$, 
\begin{equation}\label{LIMINFSPR}
s_0 = g_{0 \to k-1} (s_k) = \mbE_0\left( s_k^{Y_k} \,\mathds{1}\{Y_k<\infty\}\right) = q_0^{(k-1)} + \mbE_0\left( s_k^{Y_k} \,\mathds{1}\{0<Y_k<\infty\} \right).
\end{equation}
Suppose $\liminf_k s_k <1$. In this case there exists an infinite sequence $\{ k_i\}_{i\geq 1}$ such that $s_{k_i}<1-\varepsilon$ for all $i \geq 1$ and some $\varepsilon >0$.
For each $i \geq 1$ and $K \geq 1$, 
\[
\mbE_0\left( s_{k_i}^{Y_{k_i}} \,\mathds{1}\{0<Y_k<\infty\} \right) \leq \mbP_0(0 < Y_{k_i} < K)+ (1-\varepsilon)^K.
\]
By Lemma \ref{LemYkdic}, for any $K \geq 1$, we have $\mbP_0(0 < Y_{k_i} < K) \to 0$ as $i \to \infty$. Letting $K$ be arbitrarily large, we obtain $
\liminf_k \mbE_0( s_k^{Y_k} \mathds{1}\{ 0< Y_k < \infty  \}) = 0.$
Because $q_0^{(k)} \to q_0$ as $k \to \infty$, from \eqref{LIMINFSPR} we then obtain $s_0=q_0$. The only element $\bm{s} \in S^{[e]}$ such that $\liminf_k s_k <1$ is therefore $\bm{s}=\bm{q}$. 
\qed \medskip

Now that we have general sufficient conditions for $1-s_i \to 0$, we investigate properties of this convergence. 
The next two theorems use the following lemma.
\begin{lemma}\label{DeterminsticL}
If $\{ a_n \}_{n \geq 0}$ and $\{ b_n \}_{n \geq 0}$ are sequences of non-negative real numbers such that $a_n \in (0,1)$ for all $n \geq 0$, and $b_n \to \infty$, then 
\begin{align}\label{then1}
\limsup_n a_n^{b_n} &= \exp \{ -\liminf_n b_n(1-a_n) \},\\\label{then2}
\liminf_n a_n^{b_n} &= \exp \{ -\limsup_n b_n(1-a_n) \}.
\end{align}
\end{lemma}

\noindent\textbf{Proof:}
For any $n \geq 0$ we have 
$$
a_n^{b_n} = \left(1 - \frac{b_n(1-a_n)}{b_n} \right)^{b_n}.
$$
The result then follows from $\lim_{n \to \infty} ( 1- c/n)^n = e^{-c}$ for any $c \in \mathbb{R}$. 
\qed \medskip

The next result shows that if the entries of $\bm{q}$ converge to 1, then they converge slower than those of any other $\bm{s} \in  S\backslash\{\bm{q}\}$, whereas the entries of $\bm{\tilde q}$ converge to $1$ faster than those of any other $\bm{s} \in S\backslash\{\bm{\tilde q}, \bm{1}\}$.

\begin{theorem}\label{UniqueDecayrates}
If \eqref{LimSCond} holds then, for any $\bm{s} \in S \backslash  \{\bm{q},\bm{\tilde q},\bm{1}\}$,
\[
\lim_{k \to \infty}  \frac{1-q_i}{1-s_i}=\infty \quad \quad \mbox{ and } \quad \quad \lim_{k \to \infty} \frac{1-\tilde{q}_i}{1-s_i}=0.
\]
\end{theorem}
\noindent\textbf{Proof:}
Suppose $\bm{s} \in S \backslash \{\bm{q},\bm{\tilde q},\bm{1}\}$. In that case, by Theorem \ref{FPchar}, $\bm{q}  < \bm{s} < \bm{\tilde q}$.
 In addition, by Lemma \ref{FPequ}, for all $k\geq 0$, 
\begin{align}
s_0 &=  \mbE_0(s_k^{Y_k} \mathds{1}\{ Y_k < \infty \}) \nonumber  \\
 &= q_0^{(k-1)} + \mbE_0 \left( s_k^{Y_k} \mathds{1}\{ 0< Y_k < \infty \} \right) \label{UDE1}
 \\
&=\tilde q_0^{(k-1)} + \mbE_0\left((s_k^{Y_k}-1)\mathds{1}\{ 0<Y_k < \infty \}\right). \label{UDE2}
 \end{align}
Without loss of generality we assume that $\bm{q} < \bm{\tilde q}$, which by Corollary \ref{Link} is equivalent to $\mbP_0 \left( 0<Y_k < \infty, \, \forall \, k \geq {0}, \right)>0$.
In this case, by \eqref{UDE1} and the fact that $q_0^{(k)} \to q_0$, we have 
\begin{align*}
s_0=q_0 \quad &\Leftrightarrow \quad \lim_{k \to \infty} \mbE_0\left. \left( s_k^{Y_k}\, \right| \, 0< Y_k < \infty  \right)=0.
\end{align*}
Because $s_k^{Y_k}$ is nonnegative and uniformly bounded by 1, we can write
\begin{align*}
s_0=q_0 \quad &\Leftrightarrow \quad \mbP_0\left. \left( s_k^{Y_k} \to 0\, \right| \, \forall \, k \geq {0}, \, 0< Y_k < \infty  \right)=1.
\end{align*}
By Lemma \ref{LemYkdic} we may then apply Lemma \ref{DeterminsticL} to obtain
\begin{align*}
s_0=q_0 \quad &\Leftrightarrow \quad \mbP_0\left. \left( {Y_k}(1-s_k) \to \infty \, \right| \, \forall \, k \geq {0}, \, 0< Y_k < \infty  \right)=1,
\end{align*}hence $\lim_{k \to \infty}  (1-q_i)/(1-s_i)=\infty$.
Using \eqref{UDE2}, a similar argument yields
\begin{align*}
s_0=\tilde{q}_0 \quad  &\Leftrightarrow \quad  \mbP_0\left. \left({Y_k}(1-s_k) \to 0 \, \right| \, \forall \, k \geq {0}, \, 0<Y_k < \infty \right)=1,
\end{align*} and $\lim_{k \to \infty} (1-\tilde{q}_i)/(1-s_i)=0$.
\qed \medskip

The next theorem demonstrates that the rate at which $1-s_i$ decays is closely linked to the asymptotic growth of $\{Y_k \}$. {In this context,} we define a \emph{growth rate} to be a sequence of real numbers $\{ C_k \}_{k \geq 0}$ such that 
$$
\lim_{k \to \infty} \frac{Y_k}{C_k} = W ( \{C_k\} ) \, \, \, \, \text{exists a.s.}, $$where $W ( \{C_k\} )$ is a non-negative, potentially defective, random variable with $ \mbP(0 < W( \{ C_k \} ) < \infty) >0.$ 
We let 
\[
g_{W(\{C_k \})} (z) = \mbE_0\left( z^{W (\{ C_k \})} \mathds{1} \left\{ W (\{ C_k \})< \infty \right\} \right).
\]
Growth rates of non-defective GWPVEs ($\bm{\tilde{q}}=\bm{1}$) have been studied by a number of authors. Although it is natural to assume that $\{\mbE_0[Y_{k+1}]\}_{k\geq 0}$
 is a growth rate, it may not always be the case. Sufficient conditions for $\{\mbE_0[Y_{k+1}]\}$ to be a growth rate are given in \cite{Jag74}, and conditions for it to be the only distinct growth rate are discussed in \cite{DSo94,Big92,Kersting2017}. Examples of GWPVEs with multiple growth rates can be found in \cite{Fos76,Mac83}. 

\begin{theorem}\label{RateDecay}
Suppose \eqref{LimSCond} holds. If $\vc s\in S \backslash \{  \bm{1} \}$ and there exists some growth rate $\{ C_k \}$ such that $g_{W(\{C_k \})} (0)< s_0 < g_{W(\{C_k \})} (1)$, then
$$\lim_{k \to \infty} (1-s_k) C_k = c \in ( 0, \infty),$$ 
where $c$ is such that $s_0= g_{W(\{C_k \})} ( e^{-c})$.
\end{theorem}

\noindent\textbf{Proof:}
By the arguments in the proof of Lemma \ref{FPequ}, any $\bm{s} \in S \backslash \{  \bm{1} \}$ is such that $s_i < 1$ for all $i \in \mathcal{X}$, and $\bm{s} \in S^{[e]}$. Therefore, for all $k\geq 1$,
\begin{equation}\label{GRPeq1}
s_0 = \mbE_0\left( s_k^{Y_k} \mathds{1}\{ Y_k < \infty \} \right) = \mbE_0 \left( s_k^{Y_k} \right),
\end{equation}
which can be rewritten as
\small{\begin{equation}\label{GRPeq2}
s_0 = \mbE_0 \left( \left(s_k^{C_k}\right)^{Y_k/C_k} \mathds{1} \{W(\{C_k \}) = 0 \} \right) + \mbE_0 \left( \left(s_k^{C_k}\right)^{Y_k/C_k} \mathds{1} \{ 0< W(\{C_k \}) < \infty \}  \right)  \end{equation}$$
\hspace{-5cm}+ \mbE_0 \left( \left(s_k^{C_k}\right)^{Y_k/C_k} \mathds{1} \{W(\{C_k \}) = \infty \} \right). 
$$}
\normalsize
By assumption we have
\begin{equation}\label{GRPeq3}
g_{W(\{C_k \})}(0) = \mbP_0 \left( W ( \{ C_k \} ) =0 \right) < s_0 < g_{W( \{ C_k \} ) }( 1)=\mbP_0 \left( W ( \{ C_k \} ) < \infty \right).
\end{equation}
If $\liminf_k (s_k)^{C_k}=0$, then taking $\liminf_k$ in \eqref{GRPeq2} gives $
s_0 \leq \mbP_0(W(\{C_k \})=0),$
which contradicts \eqref{GRPeq3}.  A Similar argument applies to the limit superior, leading to
\begin{equation}\label{Th3C3}
{0<\liminf_k s_k^{C_k}\leq \limsup_k s_k^{C_k}<1.}
\end{equation}
By \eqref{GRPeq1} we then have
\begin{align}
s_0 &= \limsup_{k \to \infty} \mbE_0 \left( \left(s_k^{C_k}\right)^{Y_k/C_k} \right) \nonumber \\
&= \mbE_0 \left( \limsup_{k \to \infty} \left( s_k^{C_k}\right)^{Y_k/C_k} \right) \label{Th3C4}\\
&= \mbE_0 \left( \left(  \limsup_{k \to \infty} s_k^{C_k}\right)^{W(\{ C_k \} ) }\right),  \label{Th3C5}
\end{align}
where \eqref{Th3C4} follows from the dominated convergence theorem, and \eqref{Th3C5} requires \eqref{Th3C3}.
If we repeat the same argument with $\limsup$ replaced by $\liminf$, we finally obtain 
$$s_0 = g_{W(\{C_k\})}(\limsup_k s_k^{C_k}) = g_{W(\{C_k\})}(\liminf_k s_k^{C_k}).$$ 
By \eqref{LimSCond} and Theorem \ref{LimInfS} we have $s_k \to 1$ and thus through \eqref{Th3C3} we obtain $C_k \to \infty$. Lemma \ref{DeterminsticL} then gives $\lim_{k \to \infty}  (s_k)^{C_k} = e^{-c}$, where $c = \lim_{k \to \infty} (1-s_k) C_k$. This means
$s_0 = \mbE_0 \left( e^{-c W(\{C_k \}) }  \right) = g_{W( \{ C_k \} ) } (e^{-c})$.
\qed \medskip

We conclude this section with a summary of our findings on the set $S$.
The set $S$ is made up of a continuum of elements whose minimum is $\bm{q}$ and whose maximum is $\bm{\tilde q}$, with the additional fixed point $\bm{1}$. 
Under Condition \eqref{LimSCond}, for any $\bm{s} \in S$, with the possible exception of $\bm{q}$, we have $1-s_i \to 0$ as $i \to \infty$. 
The decay rates of $1-q_i$ and $1-\tilde q_i$ are unique, whereas the intermediate elements $\bm{q}< \bm{s} < \bm{\tilde q}$ may share one or several decay rates, which have a one-to-one correspondence with the growth rates of $\{ Y_k \}$. 
Furthermore, these intermediate elements completely specify the generating functions $g_{W(\{ C_k \})}(\cdot)$ and thereby the distributions of $W( \{ C_k \})$.
This gives a physical meaning to the intermediate elements: in short, they describe the evolution of $\{ Y_k \}$ when there is partial extinction without global extinction. 
While this physical interpretation is in terms of the growth of $\{ Y_k \}$, we expect that it is closely related to the growth of $\{ | \bm{Z}_n | \}$.

\section{Extinction Criteria} \label{SecEC}

While there exist several well-established partial extinction criteria, determining a global extinction criterion when $\bm{\tilde{q}}=\bm{1}$ remains an open question. 
When $\bm{\tilde q}=\bm{1}$, the embedded GWPVE $\{Y_k\}$ is non-explosive, and we can directly apply known extinction criteria for GWPVEs. 
These criteria are generally expressed in terms of the first and second factorial moments
$$
\mu_k := g'_k(1) \quad  \mbox{ and } \quad a_k := g''_k(1),\qquad k\geq 0.
$$
The next lemma provides recursive expressions for these moments in terms of those of the offspring distributions of $\{ \bm{Z}_n \}$.
We let 
\begin{equation}\label{mik}
m_{i\to k}:=\mbE_i[Y_{k+1}]=g'_{i\to k}(1)=\prod_{j=i}^{k} \mu_j,
\end{equation}
\[
{G}'_{k,i}(\bm{s}) := \left. \frac{\partial {G_k}( \bm{u} )}{\partial u_i} \right|_{\bm{u}=\bm{s}}, \; {G}''_{k,{ij}}(\bm{s}) := \left. \frac{\partial^2 {G_k}( \bm{u} )}{\partial u_i \partial u_j} \right|_{\bm{u}=\bm{s}},\; {A_{k,ij}:=G''_{k,ij}(\bm{1})},
\]and we take the convention that $\prod^{k-1}_{i=k} {\cdot}=1$ and $g_{k+1 \to k}(s)=s$.

\begin{lemma}\label{Embedded mean}
Suppose $\bm{\tilde{q}}=\bm{1}$, then  
\begin{equation}\label{mu_in}
\mu_{0} =\frac{M_{0,1}}{1-M_{0,0}} \quad \mbox{ and } \quad 
a_{0} = \frac{\mu_{0}^2 A_{0,00}+ A_{0,11} +2 \mu_{0}A_{0,01}}{1-M_{0,0}},
\end{equation}
and for $k \geq {1}$,
\begin{equation}\label{firstmu}
\mu_k=\frac{M_{k,k+1}}{1-\sum_{i={0}}^k M_{k,i}\,  m_{i\to k-1} },
\end{equation}
and 
\begin{equation}\label{seconda}
a_k=  \frac{\sum_{i={0}}^k M_{k,i} \sum_{j=i}^{k-1} a_j \, m_{i\to j-1}\left( \prod_{\ell=j+1}^{k} \mu^2_\ell \right)+\sum_{i={0}}^{k+1} \sum_{j={0}}^{k+1}\, m_{i\to k}\,m_{j\to k} \,A_{k,ij}}{1-\sum_{i={0}}^k M_{k,i} \, m_{i\to k-1}}.
\end{equation}
\end{lemma}
\noindent\textbf{Proof:}
By Lemma \ref{Genlink}, {for any $k\geq 0$},
\begin{align}\nonumber
g'_k(s) &= \frac{d}{ds} \left[ G_k( g_{{0} \to k} (s), \dots , g_{k+1 \to k}(s)) \right] \\\label{first_der}
&= \sum_{i={0}}^{k+1} g'_{i \to k} (s) G'_{k,i}(g_{{0}  \to k} (s), \dots , g_{k+1 \to k}(s)),
\end{align}
where 
$g'_{i \to k}(s)= \prod_{j=i}^k g'_j(g_{j+1 \to k}(s)).$
The assumption $\bm{\tilde{q}}=\bm{1}$ implies $g_{i \to k} (1)=1$ {for all $i,k$,} and therefore
$
\mu_k = {g'_k(1)}= \sum_{i={0}}^k M_{k,i} \,m_{i\to k} +M_{k,k+1},
$
which leads to {the expression for $\mu_0$ and the recursive }Equation \eqref{firstmu}.

Next, by differentiating \eqref{first_der} with respect to $s$, we obtain
\begin{align*}
g''_k(s) &= \sum^{k+1}_{i={0}} g''_{i \to k}(s)G'_{k,i}(g_{0 \to k} (s), \dots , g_{k+1 \to k}(s)) \\
&\quad + \sum^{k+1}_{i={0}} g'_{i \to k}(s) \sum_{j={0}}^{k+1} g'_{j \to k}(s) G''_{k,ij}(g_{0 \to k} (s), \dots , g_{k+1 \to k}(s)),
\end{align*}
where, {for $0\leq i\leq k$},
\[
g''_{i \to k}(s) =\sum^{k}_{j=i} \left( { \prod^{j-1}_{\ell=i  }} g'_{\ell}(g_{\ell+1 \to k}(s)) \right)g''_j(g_{j+1 \to k}(s)) \left( \prod^k_{\ell=j+1 } g'_{\ell}(g_{\ell+1 \to k}(s)) \right)^{2}.
\]
This implies
\small{$$\small{
a_k ={g''_k(1)=} \sum_{i={0}}^k M_{k,i} \sum_{j=i}^k a_j \,m_{i\to j-1}\left( \prod_{\ell=j+1}^{k} \mu^2_\ell \right)+ \sum_{i={0}}^{k+1} \sum_{j={0}}^{k+1} \,m_{i\to k}\,m_{j\to k}\, A_{k,ij},}$$}\normalsize
which gives,
\begin{align*}
a_k &\left( 1-\sum_{i={0}}^k M_{k,i} \,m_{i\to k-1}\right)= \\
&\sum_{i={0}}^k M_{k,i} \sum_{j=i}^{k-1} a_j \,m_{i\to j-1}\left( \prod_{\ell=j+1}^{k} \mu^2_\ell \right)+ \sum_{i={0}}^{k+1} \sum_{j={0}}^{k+1} \,m_{i\to k}\,m_{j\to k} \,A_{k,ij},
\end{align*}
leading to the expression for $a_0$ and the recursive Equation \eqref{seconda}. 
\qed \medskip 

When $\bm{\tilde{ q}}=\vc 1$ is not assumed, the recursive expressions \eqref{mu_in}--\eqref{seconda} can still be used to compute two sequences, which may not correspond to the first and second factorial moments of the progeny distributions of $\{ Y_k\}$, but which we shall even so denote by $\{\mu_k \}$ and $\{a_k\}$.
For these sequences to correspond to well defined moments, their elements must be non-negative and finite, that is, the denominator common to \eqref{firstmu} and \eqref{seconda} must be strictly greater than $0$ for all $k \geq 0$. Thus, if we let
\begin{equation*}
x_k : = \sum_{i={0}}^{k} M_{k,i} \,m_{i\to k-1},  
\end{equation*}
we require 
\begin{equation}\label{Algpc}
 0 \leq x_k <1\qquad \textrm{for all $k \geq 0$.}
\end{equation}
 By giving a physical interpretation to $x_k$, we now show that, in the irreducible case, \eqref{Algpc} holds if and only if $\bm{\tilde{q}}=\bm{1}$. Note that, if there exists $k$ such that $x_k=1$, then $\bm{\tilde{q}}<\bm{1}$, as justified in the proof of the lemma.

\begin{lemma}\label{Algorithmic Partial}
If $\{ \bm{Z}_n \}$ is irreducible then $\bm{\tilde{q}}=\bm{1}$ if and only if $0\leq x_k <1$ for all $k \geq 0.$
\end{lemma}
\noindent\textbf{Proof:}
For any $k\geq 0$ we embed a process 
$\{E_n^{(k)}( \bm{\tilde{Z}}^{(k)}_n) \}$ in $\{ \bm{\tilde{Z}}_n^{(k)} : \varphi_0=k\}$ by taking all type-$k$ individuals that appear in $\{ \bm{\tilde{Z}}_n^{(k)} \}$ and defining the direct descendants of these individuals as the closest (in generation) type-$k$ descendants in $\{ \bm{\tilde{Z}}^{(k)}_n \}$; the process $\{E_n^{(k)}( \bm{\tilde{Z}}^{(k)}_n) \}$ evolves as a single-type Galton-Watson process that becomes extinct if and only if type $k$ becomes extinct in $\{ \bm{\tilde{Z}}^{(k)}_n \}$.
Because $M_{i,i+1}>0$ for all $i \geq 0$, the extinction of type $k$ {in $\{ \bm{\tilde{Z}}^{(k)}_n\}$} is almost surely equivalent to the extinction of the whole process $\{ \bm{\tilde{Z}}^{(k)}_n\}$. Hence, for any $k \geq 0$,
\[
\bm{\tilde{q}}^{(k)} < \bm{1} \quad \text{ if and only if } \quad m_{E_n^{(k)}( \bm{\tilde{Z}}^{(k)}_n)}>1,
\]
where $m_{E_n^{(k)}( \bm{\tilde{Z}}^{(k)}_n)}$ is the mean number of offspring born to an individual in $\{ E_n^{(k)}(\bm{\tilde Z}^{(k)}) \}$.
The value of $m_{E_n^{(k)}( \bm{\tilde{Z}}^{(k)}_n)}$ is obtained by taking the weighted sum of all first return paths to $k$ in the mean progeny representation graph of $\{ \bm{\tilde{Z}}^{(k)}_n \}$. 
By conditioning on the progeny of an individual of type $k$ in $\{ \bm{\tilde Z}^{(k)}_n\}$, the lower-Hessenberg structure then leads to
\[
m_{E_n^{(k)}( \bm{\tilde{Z}}^{(k)}_n)} =  M_{k,0} \,m_{0\to k-1}+ M_{k,1} \,m_{1\to k-1}+ \dots + M_{k,k}=x_k.
\]
Thus, if $0<x_k < 1$ for all $k \geq 0$ then $\bm{\tilde{q}}^{(k)}=\bm{1}$ for all $k$, and therefore by Theorem \ref{PartConv},  $\bm{\tilde{q}}=\lim_{k \to \infty} \bm{\tilde{q}}^{(k)}=\bm{1}$. 
Similarly, if there exists $k$ such that $x_k>1$, then $\bm{\tilde{q}}\leq \bm{\tilde{q}}^{(k)}<\bm{1}$. 
Now suppose there exists $k$ such that $x_{k}=1$. Then by the irreducibility of $\{ \bm{Z}_n \}$ there exists $k^*>k$ such that there is a first return path with strictly positive weight of the form $k\rightarrow k+1\rightarrow \cdots\rightarrow k^*\rightarrow\cdots \rightarrow k$ in the mean progeny representation graph of $\{ \bm{\tilde{Z}}^{(k^*)}_n \}$. This implies
\[
m_{E_n^{(k)}( \bm{\tilde{Z}}^{(k^*)}_n)}>m_{E_n^{(k)}( \bm{\tilde{Z}}^{(k)}_n)}=1,
\]
and hence $
\bm{\tilde{q}} \leq \bm{\tilde{q}}^{(k^*)}<\bm{1}.$
\qed \medskip

Combining Lemmas \ref{Embedded mean} and \ref{Algorithmic Partial} with \cite[Theorem 1]{Kersting2017}, which to the authors' knowledge is the most general extinction criterion currently available for GWPVEs, we obtain for an irreducible LHBP:

\begin{theorem}\label{KerstingExtCrit}
If $\{ \mu_k \}$ and $\{ a_k \}$ are given by \eqref{mu_in}--\eqref{seconda}, then 
\begin{equation}\label{PEC}
0<\mu_k < \infty \quad \forall \, k \geq 0 \quad \Leftrightarrow \quad \bm{\tilde q}=\bm{1},
\end{equation}
and when $\bm{\tilde q}=\bm{1}$, if $\sup_{k} a_k/\mu_k < \infty$ and $\inf_k  \sum_{\{\bm{v}: v_{k+1} \geq 2\}} p_k(\bm{v}) >0$, then
\begin{equation}\label{GEC}
\sum_{k=0}^\infty \dfrac{1}{m_{0\to k}} = \infty \quad \Leftrightarrow \quad \bm{q}=\bm{1},
\end{equation}where $m_{0\to k}$ is defined in \eqref{mik}.
\end{theorem}
\noindent\textbf{Proof:}
The global extinction criterion \eqref{GEC} follows from $(i)\Leftrightarrow(iv)$ in \cite[Theorem 1]{Kersting2017}. Indeed, our assumptions imply Condition (A) of that theorem, as well as $\inf_{k} a_k/\mu_k >0$.
\qed \medskip

\begin{remark}
Theorem \ref{KerstingExtCrit} demonstrates that  by computing the sequence $\{ \mu_k\}$ 
required for \eqref{GEC} we are implementing a partial extinction criterion. 
We note that it is more efficient to compute $\{ \mu_k \}$ through Lemma \ref{Embedded mean} than to evaluate the convergence norm of $M$ as the limit of the sequence of spectral radii of the north-west truncations of the mean progeny matrix $M$ (see \cite[Theorem 6.8]{Sen06}).
\end{remark}
\begin{remark}
If $\liminf_k m_{0 \to k}=0$ then, through the Markov inequality, we obtain $\bm{q}=\bm{1}$. Thus, in this case the conditions of Theorem \ref{KerstingExtCrit} do not need to be verified.
\end{remark}

When the conditions of Theorem \ref{KerstingExtCrit} do not hold, one may still be able to apply  \cite[Theorem 1]{Kersting2017} directly.  Condition (A) in that theorem holds under an assumption on the third factorial moments $g_k'''(1)$ (\cite[Condition (C)]{Kersting2017}), which can also be shown to satisfy recursive equations.
Alternatively, it may be possible to apply 
the next theorem, which corresponds to \cite[Theorem 1]{Agr75} (see for example the proof of Proposition \ref{th_ex_2}).
\begin{theorem}\label{Agresti2} If $\bm{\tilde q} = \bm{1}$ and $A_{k,ij}<\infty$ for all $k,i,j \geq 0$, then for any $1 \leq i < k$,
$$
1 - \left[ \dfrac{1}{m_{i\rightarrow (k-1)}} + \frac{1}{2} \sum^{k-1}_{j=i} \frac{g_j''(0)}{\mu_j\,m_{i\to j}}\right]^{-1} \leq q^{(k)}_i \leq 1-  \left[\dfrac{1}{m_{i\rightarrow (k-1)}} + \sum^{k-1}_{j=i} \frac{a_j}{\mu_j\,m_{i\to j}}\right]^{-1}.
$$
\end{theorem}

Roughly speaking, Theorem \ref{KerstingExtCrit} states that the boundary between almost sure global extinction and potential global survival is the expected linear growth of $\{ Y_k \}$, that is, $\mbE_0 (Y_k) = m_{0 \to k-1} = Ck$, for some constant $C>0$. It is however not immediately clear how to interpret this criteria in terms of the expected growth of the original LHBP $\{ \bm{Z}_n \}$.
The next theorem develops a link between the expected growth of $\{ Y_k \}$ and the exponential growth rate of the mean total population size in $\{\vc Z_n\}$,
\begin{equation}\label{xiM}
\xi(M):=\liminf_n \sqrt[n]{\mbE_i | \bm{Z}_n |} = \liminf_n \sqrt[n]{(M^n\vc 1)_{i}},
\end{equation} which, when $M$ is irreducible, is independent of $i$.
We note that in an irreducible MGWBP with finitely many types $\xi(M)=\nu(M)$, whereas when there are infinitely many types it is possible that $\nu(M)<\xi(M)$.

\begin{theorem}\label{Th40}
\label{LHExpEq}Assume $\nu(M)\leq 1$. If $\xi(M)>1$, then
\begin{equation}\label{Th41}
\limsup_n { \sqrt[n]{m_{0\to n}}} \geq \xi(M),
\end{equation}
and if $\xi(M)<1$, then
\begin{equation}\label{Th43}
\liminf_n {\sqrt[n]{m_{0\to n}}} \leq \limsup_n \left( \mbE_0 |\bm{Z}_n |  \right)^{1/n}.
\end{equation}
\end{theorem}
\noindent\textbf{Proof:}
We have ${ m_{0\to( n-1)}=\sum_{k=0}^n  \mbE_0({Z}_{n,k}) \, m_{k\to( n-1)}}$, where $m_{n\to( n-1)}:=1$,
which gives
\begin{equation}\label{Th4p1}
\mbE_0 | \bm{Z}_n | \inf_{0 \leq k \leq n} {m_{k\to( n-1)}} \leq {m_{0\to( n-1)}}\leq \mbE_0 | \bm{Z}_n| \sup_{0\leq k \leq n} {m_{k\to( n-1)}}.
\end{equation}
Now suppose $\xi(M)>1$. In order to prove \eqref{Th41} we need to show that 
\begin{equation}\label{SophieE}
\nexists \, n_0<\infty \, \text{ such that } \quad \inf_{0 \leq k \leq n} {m_{k\to( n-1)}}< 1 \quad \forall \,  n>n_0.
\end{equation}
Indeed, if \eqref{SophieE} holds, because $m_{n\to( n-1)}:=1$ we have\\ $\limsup_n \inf_{0 \leq k \leq n} {m_{k\to( n-1)}}=1,$
and thus {by \eqref{Th4p1},}
\begin{align*}
\limsup_n {\sqrt[n]{m_{ 0\to( n-1)}}} &\geq \limsup_n \left( \mbE_0 | \bm{Z}_n | \inf_{0 \leq k \leq n} { m_{k\to( n-1)}} \right)^{1/n} \\
&\geq \liminf_n \left( \mbE_0 | \bm{Z}_n | \right)^{1/n} \limsup_n \left( \inf_{0 \leq k \leq n} {m_{k\to( n-1)}} \right)^{1/n} \\
&= \xi(M).
\end{align*}
To show \eqref{SophieE} assume there exists
$n_0 := \sup \left\{ n: \inf_{0 \leq k \leq n} {m_{k\to( n-1)}}= 1 \right\}<\infty,$ 
and observe that for any $n\geq 0$ the recursion 
\[
\inf_{0 \leq k \leq n} {m_{k\to( n-1)}}= \min \left\{ \left(\inf_{0 \leq k \leq n-1} {m_{k\to( n-2)}} \right) \mu_{n-1}, \,\,\mu_{n-1} ,\, \, 1 \right\}
\]
holds.
This implies that for all $n>n_0$, $$\inf_{0 \leq k \leq n} {m_{k\to( n-1)}}= \left(\inf_{0 \leq k \leq n-1} {m_{k\to( n-2)}} \right) \mu_{n-1},
$$
and
$\inf_{0 \leq k \leq n}{m_{k\to( n-1)}} = {m_{n_0\to( n-1)}},$
which gives
\[
{{m_{0\to n}}}{\left(\inf_{0 \leq k \leq n} {m_{k\to n}}\right)^{-1} }={m_{0\to( n_0-1)}}, \quad \text{ for all } n>n_0.
\]
By Equation \eqref{Th4p1} we then have $
\mbE_0 | \bm{Z}_n | \leq {m_{0\to( n_0-1)}},$  for all $n > n_0$,
which contradicts the fact that $\xi(M)>1$ and shows \eqref{SophieE}. When $\xi(M)<1$ a similar argument can be used to obtain \eqref{Th43}.
\qed \medskip

By Theorem \ref{Th40}, if both $\lim_n { \sqrt[n]{m_{0\to n}}}$ and $\lim_n \left( \mbE_0 |\bm{Z}_n |  \right)^{1/n}$ exist (which is the case in our illustrative examples), then by the \textit{root test for convergence},
$$
\xi(M)>1 \Rightarrow   \sum^\infty_{j=0} \dfrac{1}{m_{0\to j}} < \infty\quad\text{and}\quad\xi(M)<1 \Rightarrow   \sum^\infty_{j=0} \dfrac{1}{m_{0\to j}} = \infty.$$
Thus, if $\xi(M) \neq 1$ then in Theorem \ref{KerstingExtCrit}
$\sum^\infty_{j=0} 1/m_{0\to j} = \infty$
 may be replaced by $\xi(M)<1$.
One contribution of Theorem \ref{KerstingExtCrit}, which is motivated by the examples in \cite{Zuc09},
is to provide an extinction criterion applicable even when $\xi(M)=1$, as we demonstrate in Example 2.

\section{Illustrative Examples}
We now illustrate the results of the previous section through two examples.
Example 1 demonstrates that the mean progeny matrix $M$ is not sufficient to determine whether $\vc q<\vc 1$ or $\vc q=\vc 1$. 
This fact was highlighted in \cite[Example 4.4]{Zuc11}, however, in that example, the process behaves asymptotically as a GWPVE because $\sum_{j \neq i+1} M_{i,j} \to 0$ as $i \to \infty$. In addition, the proof relies on an explicit expression of the progeny generating vector.
Through Example 1 we provide a streamlined proof which applies to a significantly broader class of branching processes. 

In Example 2 we apply Theorem \ref{KerstingExtCrit} to a LHBP with $\xi(M)=1$. This example also motivates Section 7 on strong and non-strong local survival.

The proofs related to the examples are collected in Appendix B.

\smallskip
\noindent
\textbf{Example 1. }\label{PSEx6}
Consider a LHBP $\{\bm{Z}_n \}$ with mean progeny matrix
\begin{equation}\label{Tridiagonal mpg}
M=\begin{bmatrix}
b & c & 0   & 0   &    0&  \dots  & \\
a & b & c & 0   & 0   &    & \\
 0  & a & b & c & 0   &    & \\
 0  &   0 & a & b & c &    & \\
 \vdots  &    &    & \ddots & \ddots &\ddots  
\end{bmatrix},
\end{equation}
and progeny generating vector $\bm{G}( \cdot)$. We assume that $a,c>0$ and that there exists a constant $B<\infty$ such that \begin{equation}\label{assex1} A_{k,ij}=\left. \frac{\partial^2{G_k}( \bm{s} )}{\partial s_i \partial s_j} \right|_{\bm{s}=\bm{1}}
 \leq B\quad\textrm{ for all $k,i,j\geq 0$.}\end{equation}
Apart from these assumptions, we impose no other condition on $\{\bm{Z}_n \}$. 
We now consider a modification of $\{ \bm{Z}_n \}$, which we denote by $\{\bm{Z}^{\langle u \rangle}_n \}$ for some parameter $u \geq 1$, whose progeny generating vector, $\bm{G}^{\langle u \rangle}(\bm{s})$, is given by
\begin{equation}\label{Modified GenF}
G^{\langle u \rangle }_i(s_{i-1},s_i,s_{i+1}) = \frac{1}{\lceil u^i \rceil } G_i(s_{i-1},s_i,s_{i+1}^{\lceil u^i \rceil })+ \left(1-\frac{1}{\lceil u^i \rceil }\right)G_i(s_{i-1},s_i,1), \; i \geq 0.
\end{equation}
 This modification decreases the probability that a type-$i$ individual has any type-$(i+1)$ offspring by a factor of $1/\lceil u^i \rceil$, but when the type-$i$ individual does have type-$(i+1)$ offspring, their number is increased by a factor of $\lceil u^i \rceil$, which causes the mean progeny matrix to remain unchanged. Before providing results on the extinction of $\{\bm{Z}^{\langle u \rangle}_n \}$ we require the following lemma on branching processes with the tridiagonal mean progeny matrix \eqref{Tridiagonal mpg}.
\begin{lemma}\label{Try Lem}
Suppose $\{ \bm{Z}_n \}$ has a mean progeny matrix given by \eqref{Tridiagonal mpg}, then $\bm{\tilde{q}}=\bm{1}$
 if and only if 
\begin{equation}\label{patri}
{b<1} \quad \mbox{and} \quad (1-b)^2-4ac \geq 0,
\end{equation}  and when \eqref{patri} holds,
\begin{equation}\label{lim_mu}
\mu_k \nearrow \mu:= \frac{1-b - \sqrt{(1-b)^2-4ac}}{2a}\quad \text{as $k \to \infty$}.
\end{equation}
\end{lemma}
Note that the partial extinction criterion \eqref{patri} was given previously in \cite{haut12} and is implied by \cite[Theorem 1]{Big91}.
We are now in a position to {characterise the global} extinction probability of $\{\bm{Z}^{\langle u \rangle}_n \}$. 
\begin{proposition}\label{th_ex_2}
Consider the branching processes $\{ \bm{Z}^{\langle u \rangle}_n \}$ defined in Example 1, and suppose $b<1$ and $(1-b)^2-4ac>0$.
If $\mu < 1$ then $\bm{q}=\bm{1}$, whereas if $\mu  \geq 1$, then 
\[
u>\mu \,\,\, \Rightarrow \,\,\, \bm{q}=\bm{1} \quad \quad \mbox{and}  \quad\quad u<\mu \,\,\, \Rightarrow \,\,\, \bm{q}<\bm{1}, 
\]
where $\mu$ is given in \eqref{lim_mu}.
\end{proposition}

 An important sub-case of Example 1 is $u=1$, the set of unmodified branching processes.
Note that this is the only case where the second moments of the offspring distributions are uniformly bounded.
  For this subclass of processes, when combined with Lemma \ref{LemYkdic}, Proposition \ref{th_ex_2} yields 

\begin{corollary}\label{uequal1}
If $u=1$ and \eqref{LimSCond} holds then $\bm{q}=\bm{1}$ if and only if $\mu \leq 1$. 
\end{corollary}

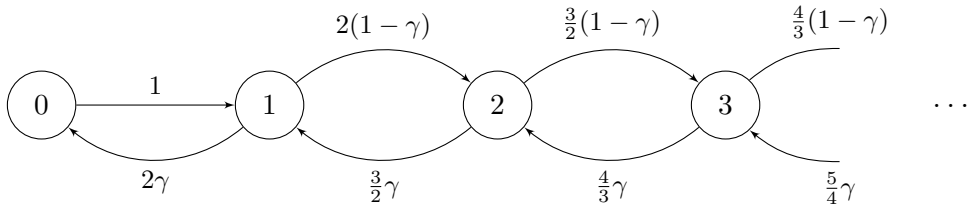
\begin{figure}
\begin{tikzpicture}

\tikzset{vertex/.style = {shape=circle,draw,minimum size=1.3em}}
\tikzset{edge/.style = {->,> = latex'}}


\node[vertex, minimum size=.9cm] (0) at  (0,0) {\small $0$};
\node[vertex, minimum size=.9cm] (1) at  (3,0) {\small $1$};
\node[vertex, minimum size=.9cm] (2) at  (6,0) {\small $2$};
\node[vertex, minimum size=.9cm] (3) at  (9,0) {\small $3$};

\node at (12,0){$\dots$};



%

\draw[edge,above] (0) to[bend left=0] node {\footnotesize $1$ } (1);

\draw[edge,above] (1) to[bend left=40] node {\footnotesize $2(1-\gamma) $ } (2);
\draw[edge,above] (2) to[bend left=40] node {\footnotesize $\frac{3}{2}(1-\gamma) $ } (3);

\draw[edge,below] (1) to[bend left=40] node {\footnotesize $2\gamma$ } (0);
\draw[edge,below] (2) to[bend left=40] node {\footnotesize $\frac{3}{2}\gamma$} (1);
\draw[edge,below] (3) to[bend left=40] node {\footnotesize $\frac{4}{3}\gamma$} (2);


\draw[above] (3) to[bend left=20, pos=1] node {\footnotesize $\frac{4}{3}(1-\gamma)$} (10.5,0.75);
\draw[edge,below] (10.5,-0.75) to[bend left=20,pos=0] node {\footnotesize $\frac{5}{4}\gamma$}  (3);

%

\end{tikzpicture}
\caption{The mean progeny representation graph corresponding to Example 2.}
\label{LHBPex1}
\end{figure}

 \noindent
\textbf{Example 2.}\label{PSEx7}
Let $\{ \bm{Z}_n \}$ have a mean progeny matrix $M$ such that $M_{0,1}=1$, and for $i \geq 1$,
\begin{equation}\label{Ex3LHBMP}
M_{i,i-1} = \gamma \frac{i+1}{i} \quad \mbox{ and } \quad M_{i,i+1} = (1-\gamma) \frac{i+1}{i}, \quad 0 \leq \gamma \leq 1,
\end{equation}
with all remaining entries being 0. The mean progeny representation graph corresponding to this process is illustrated in Figure \ref{LHBPex1}. We assume that there exists $B<\infty$ such that $A_{k,ij}  \leq B$ for all $i,j,k \geq 0$ and that $\inf_k \sum_{\bm{v}: v_{k+1}\geq 2}p_k(\bm{v})>0$.

For this example, it is not difficult to show that $\xi(M)=1$ if and only if $\nu(M)\leq 1$, which is the case for a range of values of $\gamma$, as we shall see.

\begin{proposition}\label{LHEx3s}
For the set of branching processes described in Example~2, $\bm{q}=\bm{1}$ if and only if $\gamma=0$.
\end{proposition}
Proposition \ref{LHEx3s} states that the process experiences almost sure global extinction if and only if 
type-$i$ individuals can only have type-$(i+1)$ offspring, that is, if it coincides exactly with the embedded GWPVE.

\begin{figure}
\centering
\includegraphics[width=10cm]{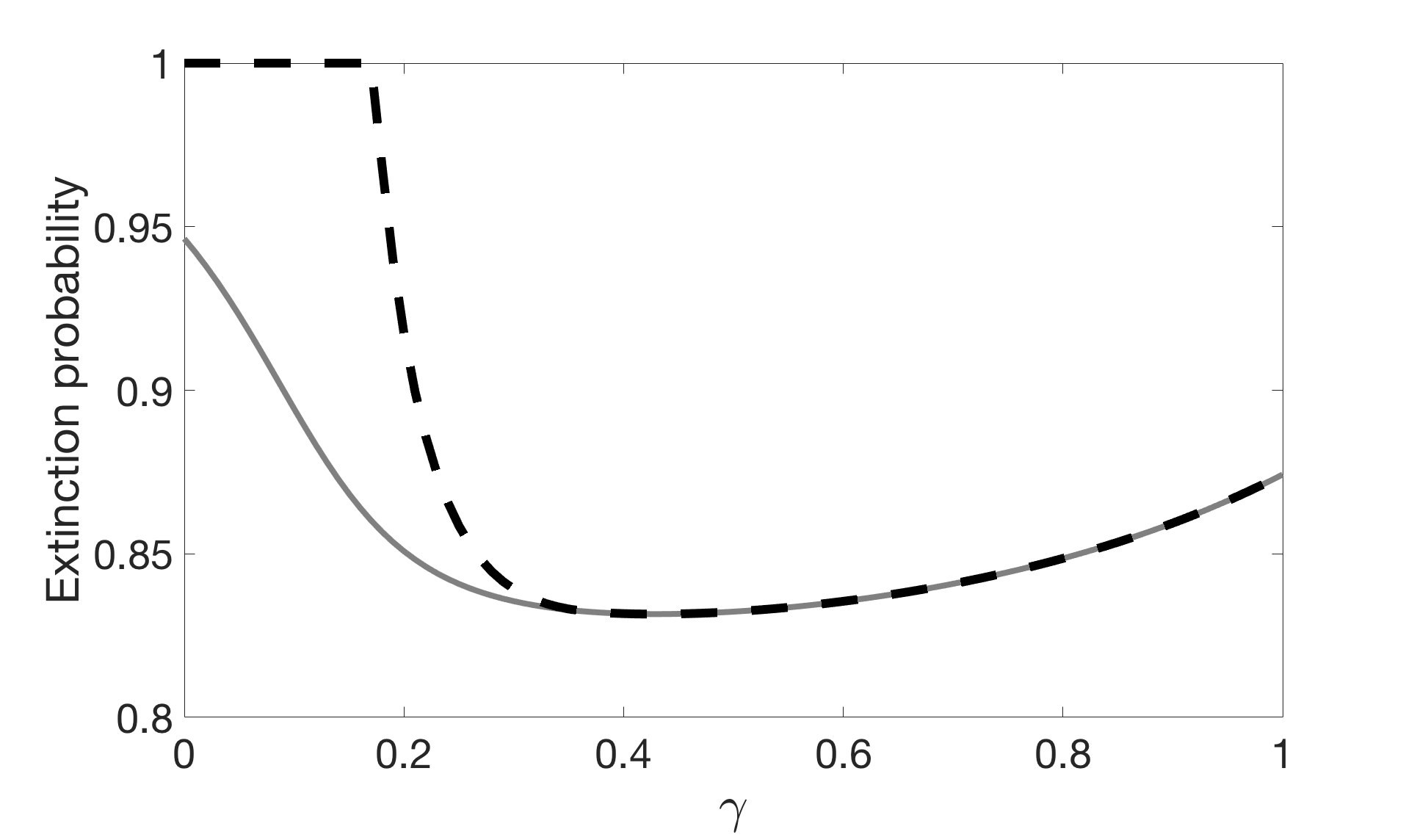}
\caption{The extinction probabilities $q_0^{(8000)}$ (solid) and $\tilde{q}_0^{(8000)}$ (dashed) for $\gamma \in [0,1]$.}
\label{LHEx3Plot}
\end{figure}

We choose
\[
G_k(\bm{s})=
\begin{cases}
\frac{1}{4}s_1^4 + \frac{1}{4}, \quad & k=0 \\
\frac{k+1}{4k}(\gamma s_{k-1} + (1-\gamma)s_{k+1})^4 + \frac{3k-1}{4k}, \quad & k\geq 1,
\end{cases}
\]
which satisfies \eqref{Ex3LHBMP}, and in Figure \ref{LHEx3Plot} we plot $q_0^{(8000)} \approx q_0$ and $\tilde{q}_0^{(8000)} \approx \tilde{q}_0$ for $\gamma \in [0,1]$. Although we proved that $q_0=1$ when $\gamma=0$, we observe that $q_0^{(8000)} \approx 0.95$ for this value of $\gamma$. This is because, when $\gamma=0$, Theorem \ref{Agresti2} implies
\[
q_0 - q_0^{(k)} \sim \left( \sum_{\ell=0}^k \frac{1}{\ell} \right)^{-1} \sim \log^{-1}(k),
\]so the convergence of $\vc q^{(k)}$ to $\vc q=\vc 1$ is slow. For GWPVEs with $\vc q<\vc 1$, little attention has been paid to this convergence rate in the literature, so for this example not much can be said when $\gamma>0$.
Using Lemmas \ref{Embedded mean} and \ref{Algorithmic Partial} we numerically determine that $\bm{\tilde{q}}=\bm{1}$ if and only if $\gamma \leq \gamma^*$ where \begin{equation}\label{gamstar}
\gamma^*= \max\{\gamma: 0<\mu_k<\infty\,\forall k\geq 0\}\approx 0.1625.\end{equation} Note that in this particular example {a sufficient condition for $\bm{\tilde{q}}=\bm{1}$ is the existence of some $k$ such that $\mu_{k} < \mu_{k-1}$} (see the proof of {Proposition} \ref{LHEx3s}). Thus, $\gamma^*$ can be {evaluated} particularly efficiently. 
Given $q^{(8000)}_0\leq q_0 \leq \tilde{q}_0 \leq \tilde{q}^{(8000)}$, by visual inspection, the curves of partial and global extinction seem to merge from some value of $\gamma$, however the cut-off is not clear and further analysis is required to pinpoint the precise value. 
   We are also interested in understanding whether this value depends only on the mean progeny matrix or whether other offspring distributions lead to different values. We address these questions in the next section.

\section{Strong local survival}\label{SecNS}

Each irreducible infinite-type branching process falls into one of the four categories $\bm{q}=\bm{\tilde{q}}=\bm{1}$, $\bm{q}<\bm{\tilde{q}}=\bm{1}$, $\bm{q}<\bm{\tilde{q}}<\bm{1}$ or $\bm{q}=\bm{\tilde{q}}<\bm{1}$. The results in the previous section deal with the classification of LHBPs with $\bm{\tilde{q}}=\bm{1}$. In the present section we build on these results to establish a method for determining whether LHBPs with $\bm{\tilde{q}}<\bm{1}$ experience strong local survival ($\vc q=\tilde{\vc q}<\vc 1$), or non-strong local survival ($\vc q<\tilde{\vc q}<\vc 1$).
Other attempts at distinguishing between these two cases can be found, for instance, in \cite{Zuc14,Zuc15} and \cite{Men97}.

For any $k\geq 0$, we partition $M$ into four components, 
$$
M=
\begin{bmatrix}
\tilde{M}^{(k)} & \bar{M}_{12} \\
\bar{M}_{21} & {}^{(k)}\tilde{M} 
\end{bmatrix},
$$
where $\tilde{M}^{(k)}$ is of dimension$(k+1) \times (k+1)$ and the other three submatrices are infinite. We then construct a LHBP branching processes on $(\Omega, \mathcal{F}, \mbP)$, denoted as $\{{}^{(k)}\bm{\tilde{Z}}_n\}$,  with mean progeny matrix ${}^{(k)}\tilde{M}$, and {global and partial} extinction probability vectors ${}^{(k)}\bm{q}$ and ${}^{(k)}\bm{\tilde{q}}$, respectively. Sample paths of $\{{}^{(k)}\bm{\tilde{Z}}_n\}$ are constructed from those of $\{ \bm{Z}_n \}$ by immediately killing all offspring of type $i \leq k$, and relabelling the types so that type $i\geq k+1$ becomes $i-k-1$. We now use $\{{}^{(k)}\bm{\tilde{Z}}_n\}$ to derive a criterion for strong local survival. In the next theorem we let $sp(\cdot)$ denote the spectral radius.

\begin{theorem}\label{Case4}
Assume that $\bm{\tilde{q}}>\bm{0}$, that there exists $k \in \mbN$ such that
\begin{equation}\label{assumptions}
sp\left( \tilde{M}^{(k)} \right)>1 \quad \mbox{and} \quad \nu \left( {}^{(k)} \tilde{M} \right) \leq 1,
\end{equation}
and that $\bar{M}_{21}$ contains a finite number of strictly positive entries.
Then there is strong local survival in $\{ \bm{Z}_n \}$ if and only if $\{{}^{(k)}\bm{\tilde{Z}}_n\}$ becomes globally extinct, that is,
\begin{equation}\label{sls}
\bm{q} = \bm{\tilde{q}} < \bm{1} \quad  \text{ if and only if }  \quad {}^{(k)}\bm{q} = {}^{(k)} \bm{\tilde{q}}= \bm{1}.
\end{equation}
\end{theorem}
\noindent\textbf{Proof:}
 We use \cite[Theorem 4.2]{Zuc15} which we restate using our notation: let $\{\vc Z_n^{(G)}\}$ and $\{\vc Z_n^{(G^*)}\}$ be two branching processes on the countable type set $\mathcal{X}$ with respective probability generating functions $\vc G(\cdot)$ and $\vc G^*(\cdot)$, and global extinction probability vectors $\vc q$ and $\vc q^*$. Let $A\subseteq\mathcal{X}$ be a non-empty subset of types and denote by $\vc q(A)$ and $\vc q^*(A)$ the respective vectors of probability of local extinction in $A$. If  $\{\vc Z_n^{(G)}\}$ and $\{\vc Z_n^{(G^*)}\}$ differ on $A$ only, that is, if $G_i(\vc s)=G_i^*(\vc s)$ for all $i\in \mathcal{X}\setminus A$ and $G_i(\vc s)\neq G_i^*(\vc s)$ for all $i\in A$, then
\begin{equation}\label{iaoi}\vc q  = \vc q(A)\quad  \Leftrightarrow \quad \vc q^* =\vc q^*(A) .\end{equation}
We apply this result with $A=\{0,1,\ldots,k\}$, $\{\vc Z_n^{(G)}\}=\{\vc Z_n\}$, and $\{\vc Z_n^{(G^*)}\}$ being such that $G_i^*(\vc s)=\vc 1$ for all $i\in A$, that is, all types in $A$ are sterile. We need to show that \eqref{iaoi} is equivalent to \eqref{sls}.

We first observe that $\vc q(A) =\vc{\tilde{q}}$ since, by \eqref{assumptions}, in $\{\vc Z_n\}$, types in $\mathcal{X}\setminus A$ are only able to survive through the presence of types in $A$. Next, since types in $A$ are sterile in $\{\vc Z_n^{(G^*)}\}$, $\vc q^* =\vc q^*(A)$ if and only if $ q_i^* =q_i^*(A)$ for all $i\geq k+1$. It is clear that $(q_{k+1}^*,q_{k+2}^*,\ldots)={}^{(k)}\bm{q}$ by construction. It remains to show that $(q_{k+1}^*(A),q_{k+2}^*(A),\ldots)=\vc 1$. We couple the process $\{{}^{(k)}\bm{\tilde{Z}}_n:\varphi_0= \ell-1\}$  and the process $\{\vc Z_n^{(G^*)}:\varphi_0=k+\ell\}$, $\ell\geq 1$.  Let $k+\bar{\ell}$ be the largest type in $\{\vc Z_n^{(G^*)}\}$ able to generate offspring in $A$. Then since $\nu ( {}^{(k)} \tilde{M} ) \leq 1$, with probability one there exists a generation $N$ such that ${}^{(k)}{\tilde{Z}}_{n, 0}+\ldots+{}^{(k)}{\tilde{Z}}_{n,\bar{\ell}}=0$ for all $n\geq N$. This implies that with probability one, $ Z_{n,k+1}^{(G^*)}+\ldots + Z_{n,k+\bar{\ell}}^{(G^*)}=0$ for all $n\geq N$, which shows $(q_{k+1}^*(A),q_{k+2}^*(A),\ldots)=\vc 1$. 
\qed \medskip 

When $\nu \left( {}^{(k)} \tilde{M} \right) \leq 1$, we may apply Theorem \ref{KerstingExtCrit} to determine whether ${}^{(k)}\bm{q}  = \bm{1}$. 
We are now in a position to answer the questions posed at the end of the previous section. 
The next result is proved in Appendix B.

\begin{proposition}\label{LHEx3final}
For the branching processes described in Example~2,
\begin{align*}
\gamma = 0 \quad &\Rightarrow \quad \bm{q}=\bm{\tilde{q}}=\bm{1} \\
\gamma \in (0,\gamma^*] \quad &\Rightarrow \quad \bm{q}<\bm{\tilde{q}}=\bm{1} \\
\gamma \in (\gamma^*,1/2) \quad &\Rightarrow \quad \bm{q}<\bm{\tilde{q}}<\bm{1} \\
\gamma \in (1/2,1] \quad &\Rightarrow \quad \bm{q}=\bm{\tilde{q}} < \bm{1},
\end{align*}where $\gamma^*$ is given in \eqref{gamstar}.
\end{proposition}

Proposition \ref{LHEx3final} demonstrates that the curves for partial and global extinction represented in Figure \ref{LHEx3Plot} merge at $\gamma=1/2$ and that this value is independent of the particular offspring distributions. At the critical value $\gamma=1/2$ there exists no $k$ satisfying \eqref{assumptions}, causing this case to remain untreated.

\section{Conclusion}

 Besides thoroughly exploring the set of fixed-points for LHBPs, we have introduced a method of classifying LHBPs 
into one of the categories $\bm{q}=\bm{\tilde{q}}=\bm{1}$, $\bm{q}<\bm{\tilde{q}}=\bm{1}$, $\bm{q}<\bm{\tilde{q}}<\bm{1}$ or $\bm{q}=\bm{\tilde{q}}<\bm{1}$.Through Examples 1 and 2 we showed that our results can be used to rigorously determine which category the process falls in; however, in practical situations where rigorous proofs may not be possible, our results can still be applied computationally as a first step in classifying the process.

The inherent assumption in LHBPs is the constraint that individuals of type $i$ cannot give birth to offspring whose type is larger than $i+m$ for $m=1$. The approach of embedding a GWPVE in the original LHBP can be extended to the case where $m$ takes any finite integer value. The resulting embedded GWPVE then becomes multitype with $m$ types. Results of Section \ref{sec_GWPVE} then naturally generalise, but those of Section \ref{fpep} rely on the characterisation of the $m$-dimensional projection sets of $S$, which is more difficult in this case. The global extinction criterion discussed in Section \ref{SecEC} would now build upon extinction criteria for multitype GWPVE, which are less developed in the literature. These questions are the topic of a subsequent paper \cite{Bra18}.

\appendix

\section*{Appendix A: Partial extinction probability}

The following result holds not only for LHBPs but in the more general setting of \cite{haut12}.

\begin{theorem}\label{PartConv}
If $\{ \bm{Z}_n \}$ is irreducible and non-singular then $\bm{\tilde q}^{(k)} \searrow \bm{\tilde q}$ componentwise as $k \to \infty$.
\end{theorem} 
\noindent\textbf{Proof:}
Fix some initial type $i \in \mathcal{X}$.
By construction, for every $n\geq 0$ and $\omega \in \Omega$, $\bm{\tilde Z}_n^{(k)}(\omega)$ is increasing in $k$, which implies $\tilde q_i^{(k)}$ is decreasing in $k$. Similarly, if $\{ \bm{\tilde Z}^{(k)}_n \}$ survives globally, then at least one type $j \in \{1, \dots, k \}$ must survive in $\{ \bm{Z}_n \}$, which implies $\tilde q_i^{(k)} \geq \tilde q_i$ for all $k$.
We may then assume $\tilde q_i <1$. 
Because $\{ \bm{Z}_n \}$ is irreducible, \cite[Corollary 1]{Bra18} implies that $\tilde q_i$ is equal to the probability that type $i$ eventually disappears from the population. 
We define a function $f^{(i)}: \mathcal{J} \to \mathcal{J}$ that takes lines of descent $(\varphi_0; i_1,j_1,y_1; \dots ; i_n, j _n, y_n)$ and deletes each triple whose type $j_{(\cdot)}$ is not equal to $i$, and we define the processes $\{ V_{\ell}^{(i)} \}_{\ell\geq 0}$ and $\{ \tilde V_{\ell}^{(i,k)} \}_{\ell\geq 0}$, whose family trees are given by $f^{(i)}(X)$ and $f^{(i)}(\tilde X^{(k)})$, respectively. These are single-type Galton-Watson processes that become extinct if and only if type $i$ becomes extinct in $\{ \bm{Z}_n \}$ and $\{ \bm{\tilde Z}^{(k)}_n \}$.
Thus, given $\varphi_0=i$ the probability that $\{ V_{\ell}^{(i)} \}$ becomes extinct is $\tilde q_i$, and the probability that $\{ \tilde V_{\ell}^{(i,k)} \}$ becomes extinct is greater than or equal to $\tilde q^{(k)}_i$ (if $\{ \bm{\tilde Z}_n^{(k)} \}$ is reducible there may be a positive chance type $i$ dies out but $\{ \bm{\tilde Z}^{(k)}_n \}$ survives globally).
%
Because $\{ \bm{Z}_n \}$ is irreducible and non-singular, $\{ V_{\ell}^{(i)} \}$ is non-singular, that is, there is positive chance that individuals in $\{ V_{\ell}^{(i)} \}$ have a total number of offspring different from 1.
Thus, with probability 1, $\{ V^{(i)}_\ell \}$ experiences extinction or unbounded growth \cite[Chapter I, Theorem 6.2]{Har63}. For any $K >0$ we then have
\begin{equation}\label{Conv1}
\lim_{\ell \to \infty} {\mbP}_i( V^{(i)}_\ell \geq K ) = 1 - \tilde{q}_i.
\end{equation}
Observe that, for any  fixed $h \in \mbN$ and $K >0$,
\small{\[
\{ \omega\in\Omega : \tilde V_h^{(i,1)}(\omega)  \geq K \} \subseteq \{ \omega\in\Omega : \tilde V_h^{(i,2)}(\omega) \geq K \} \subseteq\{ \omega\in\Omega : \tilde V_h^{(i,3)}(\omega) \geq K \} \subseteq \dots
\]}
\normalsize
and 
\begin{equation}\label{Asser2}
\lim_{k \to \infty} \{ \omega\in\Omega : \tilde V_h^{(i,k)}(\omega)  \geq K \} = \{ \omega\in\Omega : V_{h}^{(i)}(\omega) \geq K \}.
\end{equation}
To understand \eqref{Asser2} observe that if $\omega \in \{ \omega\in\Omega : V_h^{(i)}(\omega) \geq K \}$, then there exists at least $K$ lines of descent $(\varphi_0; i_1,j_1,y_1; \dots ; i_n, j _n, y_n) \in X(\omega)$ such that the type $j_n=i$ is the $h$th return to $i$ (where $n$ is not necessarily the same for each of these $K$ lines of descent).
By construction, the maximum type on each of these lines of descent is finite. Thus letting $k_m$ denote the maximum of the maximum type on $K$ arbitrarily selected such lines of descent, we see that $\omega \in  \{ \omega\in\Omega : \tilde V_h^{(i,k_m)}(\omega)  \geq K \}$.
We may now apply the monotone convergence theorem to obtain, for any $h \in \mbN$,
\begin{equation}\label{Conv2}
\lim_{ k \to \infty} {\mbP}_i( \tilde V^{(i,k)}_h \geq K ) = {\mbP}_i( V^{(i)}_h \geq K ).
\end{equation}
The probability that $\{ \tilde V_\ell^{(i,k)} \}_{\ell \geq 0}$ becomes extinct is equal to that of $\{ \tilde V^{(i,k)}_{h \ell} \}_{\ell \geq 0}$, which is less than or equal the probability of extinction $\tilde{q}_i^{(k,h,K)}$ of the Galton-Watson branching process with progeny generating function 
\[
G(s) = 1- {\mbP}_i( \tilde V^{(i,k)}_h \geq K ) + s^K{\mbP}_i( \tilde V^{(i,k)}_h \geq K ).
\]
Therefore $\tilde q^{(k)}_i \leq \tilde q_i^{(k,h,K)}$. 
Observe that for any  fixed $p \in (0,1]$ the probability of extinction in a Galton-Watson process with progeny generating function $G(s)=1-p + ps^K$ converges monotonically to $1-p$ as $K \to \infty$ (to see why note that for any $0<\eta<p$ there exists $K$ large enough to ensure $G(1-p+\eta) = (1-p) + p((1-p) + \eta)^K \leq 1-p+\eta$).
We are now in a position to show that for any $\varepsilon >0$ there exists $k_i$ such that $\tilde q^{(k_i)}_i < \tilde q_i +\varepsilon$. 
Given a process with partial extinction probability $\tilde q_i<1$ and some $0<\varepsilon<1-\tilde q_i$,
we select $K$ by setting it large enough to ensure that a Galton-Watson branching process with progeny generating function $G(s)=(\tilde q_i + \varepsilon/2) + (1-(\tilde q_i + \varepsilon/2))s^K$ has extinction probability less than $\tilde q_i + \varepsilon$.
By \eqref{Conv1}, for this value of $K$, we may select $h$ large enough to ensure $|{\mbP}_i( V^{(i)}_h \geq K ) - (1 - \tilde q_i) |<\varepsilon/4$. By \eqref{Conv2}, for these values of $K$ and $h$ we may select $k_i$ large enough to ensure $|{\mbP}_i( V^{(i)}_h \geq K )-{\mbP}_i( \tilde V^{(i,k_i)}_h \geq K )| < \varepsilon/4$. By the triangle inequality and the preceding discussion, for these values of $k_i$, $h$, and $K$, we have $\tilde q_i \leq \tilde q_i^{(k_i)} \leq \tilde q_i^{(k_i,h,K)} < \tilde q_i + \varepsilon$. The result then follows from the fact $\tilde q_i^{(k)}$ is decreasing in $k$. 
\qed
 
 \section*{Appendix B: Proofs related to the examples}

\textbf{Proof of Lemma \ref{Try Lem}:} Because \eqref{patri} holds, Lemma \ref{Embedded mean} gives 
\begin{equation}\label{L9LH}
\mu_0=\frac{c}{1-b}, \quad \mbox{ and } \quad  \mu_k =\frac{c}{1-b-a \mu_{k-1}}\quad \textrm{for all $k \geq 0$. }
\end{equation}
Because $\mu_1>\mu_0$ and
\[
\mu_k - \mu_{k-1} = \frac{c}{1-b-a \mu_{k-1}}-\frac{c}{1-b-a \mu_{k-2}},
\]
by induction the sequence $\{\mu_k\}_{k\geq 0}$ is strictly positive and increasing.  
Therefore, since $a>0$, 
$\bm{\tilde{q}}=\bm{1}$ implies that $\{\mu_k\}$ converges to a finite limit $\mu$, 
where $\mu$ satisfies the equation
$
ax^2-(1-b)x+c=0,
$
which has real solutions
\[
{x_{\pm}=} \frac{1-b \pm \sqrt{(1-b)^2-4ac}}{2a},
\]
since \eqref{patri} holds. When \eqref{patri} holds we have 
$\mu_0\leq x_{-}$ which, combined with \eqref{L9LH} and the fact that $x_{-}=c/(1-b-a x_{-}),$ implies $\mu_k\leq x_{-}$ for all $k\geq 0$, hence $\mu_k \nearrow \mu = x_{-}$. 
\qed \medskip

\noindent\textbf{Proof of Propostion \ref{th_ex_2}:} 
Let $ \Delta=(1-b)^2-4ac>0$. First, suppose $u>\mu$. In this case we have
$$
1-q^{(k)}_0 = \frac{\mbE_0(Y_k)}{\mbE_0(Y_k | Y_k >0)} \leq \frac{\mu^k}{u^{k-1} },
$$
where $\mbE_0(Y_k) \leq \mu^k$ follows from Lemma \ref{Try Lem} and $\mbE_0(Y_k | Y_k >0) \geq u^{k-1}$ follows from the fact that the minimum number of type-$k$ offspring born to a type-$(k-1)$ parent is $\lceil u^{k-1} \rceil$.
This then implies
\[
1-q_0=  1-\lim_{k \to \infty} q^{(k)}_0 \leq \lim_{k \to \infty} \frac{\mu^k}{u^{k-1}}= 0,
\]and therefore $\vc q=\vc 1$ by irreducibility.

Now suppose $1\leq  u < \mu$. Note that $A^{\langle u \rangle}_{k,ij}=A_{k,ij}$ for all $i,j$ with the exception of $A^{\langle u \rangle}_{k,(k+1)(k+1)}=\lceil u^k \rceil A_{k,(k+1)(k+1)}+ c({\lceil u^k \rceil-1})$.
Then, by Lemma~\ref{Embedded mean},
\begin{align*}
a_k &=  \frac{a \mu_k^2 a_{k-1} + \lceil u^k \rceil A_{k,(k+1)(k+1)}+ O(1)}{1- b - a \mu_{k-1}} \\
&\leq  \frac{a \mu^2 a_{k-1} + B \lceil u^k \rceil+O(1)}{1- b - a \mu} \\
&=a_{k-1} \mu \frac{a\frac{1-b - {\Delta}^{1/2}}{2a}}{1-b-a\frac{1-b - {\Delta}^{1/2}}{2a}}+ B^* u^k+O(1)\\
&=a_{k-1} \mu \frac{{1-b}- \Delta^{1/2}}{{1-b}+\Delta^{1/2}}+ B^* u^k+O(1),
\end{align*}
for all $k \geq 0$ and some $B^*<\infty$, which implies
\[
a_k = O \left( \left[ \max \left\{ u, \mu ({1-b}- \Delta^{1/2})/({1-b}+\Delta^{1/2}) \right\} \right]^k\right).
\]
By assumption, $\Delta>0$ and $u < \mu$, thus $ \max \left\{ u, \mu \frac{{1-b}- \Delta^{1/2}}{{1-b}+\Delta^{1/2}} \right\}<\mu$. Using the fact that $\mu_k \nearrow \mu$ and the \emph{root test}, we then obtain
\[
\sum^\infty_{k=0} \frac{a_k}{\mu_k \,m_{0\to k}}<\infty,
\] 
which, by the upper bound in Theorem \ref{Agresti2}, gives $q_0 <1$. 
\qed \medskip

\noindent\textbf{Proof of Corollary \ref{uequal1}:} It remains to show $\bm{q}=\bm{1}$ when $\mu=1$. Lemma~\ref{LemYkdic} implies $\mbP_0(Y_k \to 0)+ \mbP_0(Y_k \to \infty) =1$ and Lemma \ref{Try Lem} implies $\mbE_0(Y_k)= \prod^{k-1}_{i=0} \mu_i \leq 1$ for all $k$, leading to $\mbP_0(Y_k \to \infty) =0$ and the result. \qed \medskip

\smallskip 
\noindent\textbf{Proof of Proposition \ref{LHEx3s}:}
If $\gamma=0$ then $\mu_0=M_{0,1}=1$, and for $k \geq 1$, $\mu_k = (k+1)/k$. This gives $m_{0\to k} ={k+1}$, and therefore 
$$
\sum_{k= 0}^\infty \dfrac{1}{m_{0\to k}}= \sum_{k=1}^\infty \frac{1}{k}=\infty.
$$
By assumption, when $\gamma=0$ the conditions Theorem \ref{KerstingExtCrit} are satisfied, which then implies $\vc q=\vc 1$.

Now suppose $\gamma>0$. By Lemma \ref{Embedded mean}, for $k \geq 1$,
\begin{equation}\label{LHEq11}
\mu_{k}=f_k(\mu_{k-1}):= \frac{\frac{k+1}{k}(1-\gamma)}{1-\frac{k+1}{k} \gamma \mu_{k-1}}.
\end{equation}
If there exists $k$ such that $\mu_k>1/\gamma$, then by Lemma \ref{Algorithmic Partial} we have $\bm{q} \leq \bm{\tilde{q}}<\bm{1}$. Assume from now on that $\bm{\tilde{q}}=\bm{1}$, which implies that $0 \leq \mu_k \leq 1/\gamma < \infty$ for all $k$, and $\gamma<1/2$.
 Since $\mu_0=M_{0,1}=1$, using Equation \eqref{LHEq11} we can inductively show that $\mu_k \geq 1$ for all $k \geq 0$. {We then have, for any $k\geq 1$,}
\begin{equation}\label{LHExAG}
\mu_k \;\geq\; \frac{\frac{k+1}{k}(1-\gamma)}{1-\frac{k+1}{k}\gamma}  \;\geq\; 1+\frac{1}{k(1-\gamma)}. 
\end{equation}
The \emph{Raabe-Duhamel test} for convergence ensures that $
\sum_{k={0}}^\infty (1/m_{0\to k})<\infty$,
since for $k \geq 1$,
\[
k \left( \frac{(1/m_{0\to (k-1)})}{(1/m_{0\to k})} -1\right)=k(\mu_k-1) \geq \frac{1}{1-\gamma}>1.
\]

To complete the proof,  it remains to show that the condition $\sup_{k} a_k/\mu_k < \infty$ in Theorem \ref{KerstingExtCrit} holds.  
By Lemma \ref{Embedded mean}, for all $k\geq 1$,
\begin{eqnarray*}a_k &=& \frac{ a_{k-1}\,\gamma \frac{k+1}{k}  \,\mu_{k}^2 }{1-\gamma \frac{k+1}{k}  \mu_{k-1}}\\&&+\;\frac{ A_{k,(k-1)(k-1)}\mu_{k-1}^2\mu_k^2+ 2 A_{k,(k-1)(k+1)}\mu_{k-1}\mu_k+A_{k,(k+1)(k+1)}}{1-\gamma \frac{k+1}{k}  \mu_{k-1}} .\end{eqnarray*} 
Since $\bm{\tilde{q}}=\bm{1}$, the denominator is uniformly bounded away from 0; in addition, by assumption, $A_{k,ij}  \leq B<\infty$ for all $i,j,k \geq 0$, therefore
there exists some constant $K < \infty$ independent of $k$ such that
\begin{align*}
a_k \leq 
\frac{a_{k-1} \gamma \frac{k+1}{k} \mu_k^2}{1-\frac{k+1}{k} \gamma \mu_{k-1}}+K.
\end{align*}
If $\mu_k \to1$ (which we show below), then for large $k$,
\[
a_k \leq \left(\frac{\gamma}{1-\gamma} +o(1)\right) a_{k-1}+K.
\]
{Since $\gamma<1/2$}, we have $\gamma/(1-\gamma)<1$, which means that $\{a_k\}$ is a uniformly bounded sequence. Combining this with the fact that $\mu_k \geq 1$ for all $k$ implies $\sup a_k / \mu_k <\infty$, and $\vc q<\vc 1$ by Theorem \ref{KerstingExtCrit}. 

Finally, we prove that $\mu_k \to 1$. Observe that \eqref{LHEq11} implies that if  $\mu_{k} < \mu_{k-1}$ {for some $k$}, then $\mu_{k+1} < \mu_{k}$, and thus $\mu=\lim_{k \to \infty} \mu_{k}$ exists since $1{\leq} \mu_k\leq 1/\gamma$ for all $k$.  
Taking $k \to \infty$ in \eqref{LHEq11} we obtain that $\mu$ satisfies
\[
\mu = \frac{1-\gamma}{1-\gamma \mu}:=f(\mu),
\]
which means $\mu$ is either $1$ or $(1-\gamma)/\gamma>1$. 
The function $f(x)$ is convex, thus 
$f(x)>x$ for all $x>(1-\gamma)/\gamma$; in addition, by \eqref{LHEq11}, $\mu_{k+1}>f(\mu_{k})$ for all $k\geq 0$. These imply that  if $\mu_k>(1-\gamma)/\gamma$ for some $k$, then $\mu_{k+\ell}$ becomes negative for some $\ell>1$, which is a contradiction. So the sequence $\{\mu_k\}$ lives in the open interval $(1,(1-\gamma)/\gamma)$. 
Let \begin{equation*}
 v_{\pm}^{(k)}=\dfrac{1}{2\gamma}\left[\frac{k}{k+1}\pm \sqrt{\left( \frac{k}{k+1}\right)^2 - 4 \gamma (1-\gamma)}\right]
\end{equation*}be the solutions of the equation $x=f_k(x)$. By the convexity of $f_k(x)$ for all $k$, if there exists $K\geq 1$ such that $v_{-}^{(K+1)}<\mu_{K}<v_{+}^{(K+1)}$ then $\{ \mu_k \}_{k \geq K}$ is a decreasing sequence which converges to 1.  Suppose $\mu = (1-\gamma)/\gamma$. Then $\mu_{K} \geq v_+^{(K+1)}$ for some $K$. We can then construct a LHBP, $\{ \bm{Z}^*_n\}$, stochastically smaller than $\{ \bm{Z}_n \}$ by selecting a sufficiently large type $K$ and independently killing each type-$(K+1)$ child born to a type-$K$ parent with a probability carefully chosen to ensure $v_{-}^{*(K+1)}< \mu_K^*< v_{+}^{*(K+1)}$. For this modified process we have $\mu^*=1$, and repeating previous arguments, we obtain $\bm{q}<\bm{q}^*<\bm{1}$.
 \qed \medskip

\noindent\textbf{Proof of Proposition \ref{LHEx3final}:} Given Proposition \ref{LHEx3s} and Lemmas \ref{Embedded mean} and \ref{Algorithmic Partial}, it remains to show that $\bm{q}<\bm{\tilde{q}}$ for $\gamma \in (\gamma^*,1/2)$ and $\bm{q}=\bm{\tilde{q}}$ for $\gamma \in (1/2,1]$. Note that, in either case, since $\bm{\tilde{q}}<\bm{1}$, ${ \exists K_1 \; \textrm{such that $sp(\tilde{M}^{(k)})>1$\; $\forall k\geq K_1$}}.${ In addition, $$\forall x>1, \;\exists K(x)\;\textrm{s.t. $M_{k,k+1} <x(1-\gamma)$ and $M_{k,k-1}<x\gamma$ \; $\forall k \geq K(x)$}.$$ Since $\gamma \neq 1/2$, we may choose $\bar{x}>1$ small enough so that $1-4\bar{x}^2(1-\gamma)\gamma>0$. By Lemma \ref{Try Lem}, this implies that $\nu({}^{(k)}\tilde{M})<1$ for all $k\geq \bar{K}:=K(\bar{x})$, and 
\[
{}^{(\bar{K})}\mu_k \leq \frac{1-\sqrt{1-4\bar{x}^2(1-\gamma)\gamma}}{2 \bar{x} \gamma}
\]
for all $k\geq 0$, where $\{{}^{(\bar{K})}\mu_k\}_{k \geq 0}$ is computed using ${}^{(\bar{K})}\tilde{M}$. 

Assume first that $\gamma \in (1/2,1]$. Then,
$
(1-\sqrt{1-4(1-\gamma)\gamma})/(2  \gamma)<1,
$
so we may choose $x^*\leq \bar{x}$ small enough, corresponding to $K^*:=K(x^*)\geq \bar{K}$, so that ${}^{(K^*)}\mu_k<1-\varepsilon$ for all $k\geq 0$ and some $\varepsilon>0$. Hence there exists $K=\max\{K_1,K^*\}<\infty$} satisfying the conditions of Theorem \ref{Case4} with ${}^{(K)}\bm{q} = {}^{(K)} \bm{\tilde{q}}= \bm{1}$.

Now suppose $\gamma \in (c,1/2)$. If for any $K>\bar{K}$ there exists $k_1 \geq 0$ such that $
{}^{(K)}\mu_{k_1} \geq 1,$
then by the recursion \eqref{LHEq11} we have $
{}^{(K)}\mu_k \geq 1 + 1/(k(1-\gamma))$, for all $k > k_1,$
and the result is derived by repeating the steps that follow Equation \eqref{LHExAG} in the proof of  Proposition \ref{LHEx3s}. Suppose instead that there exists $K > \bar{K}$ such that ${}^{(K)} \mu_k < 1$ for all $k \geq 0$. Then by Equation \eqref{LHEq11} we have $
{}^{(K)}\mu_{k-1}< 1- 1/(\gamma(k+1)),$
which implies 
\begin{equation}\label{Th811}
\prod^\infty_{i=0} {}^{(K)} \mu_{i}=0.
\end{equation}
To show that this leads to a contradiction, we compare $M$ to a matrix $M^*$ with strictly smaller entries than $M$: $M^*$ is such that $M^*_{0,1}=1-\gamma$, and for all $k \geq1$, $M^*_{k,k-1}=\gamma$ and $M^*_{k,k+1}=1-\gamma$,
with all other entries 0. The value of $\prod^\infty_{i=0} \mu_i^*$, with $\{\mu_i^*\}_{i \geq 0}$ computed using $M^*$, then has a probabilistic interpretation: it is the probability that a simple random walk on the integers, with transition probabilities $p_{+}=1-\gamma>p_{-}=\gamma$, whose initial value is 0, never hits $-1$. When $1-\gamma>\gamma$ it is well known that this value is non-zero. By the fact that ${}^{(K)}M>M^*$ we then have ${}^{(K)}\mu_i>\mu^*_i$ for all $i \geq 0$, which implies $
\prod^\infty_{i=0} {}^{(K)}\mu_i > \prod^\infty_{i=0} \mu_i^* > 0
,$
contradicting \eqref{Th811}. \qed \medskip

\section*{Acknowledgements}
The authors are grateful to two anonymous referees for their constructive comments, which helped us to improve the manuscript. In particular, they thank the referee who pointed out the inaccuracy in \cite[Lemma 3.2]{haut12}.
The authors would also like to acknowledge the support of the Australian Research Council (ARC) through the Centre of Excellence for the Mathematical and Statistical Frontiers (ACEMS). Sophie Hautphenne would further like to thank the ARC for support through Discovery Early Career Researcher Award DE150101044.

\end{document}